\documentclass[12pt]{article}
\usepackage{amsmath}
\usepackage{amsthm}
\usepackage{amsfonts}
\usepackage{amssymb}
\usepackage{esint}
\usepackage{eucal}
\usepackage{mathrsfs}
\usepackage{graphicx,graphics}
\numberwithin{equation}{section}
\usepackage{graphicx,graphics}
\usepackage{pdfsync}
\usepackage{color}
\def \dis {\displaystyle}

\def \into {\int_\Omega}
\def \confai {-\kern -.5em\rightharpoonup}
\def \cqfd {\hfill$\Box$}

\def \de {\delta}
\def \De {\Delta}
\def \ep {\varepsilon}
\def \om {\omega}
\def \Om {\Omega}
\def \la {\lambda}

\def \ph {\varphi}

\def \NN {\mathbb N}

\def \RR {\mathbb R}

\def\C {\mathscr{C}}
\def \D {\mathscr{D}}
\def \F {\mathscr{F}}
\def \G {\mathscr{G}}

\def \L {\mathscr{L}}
\def \M {\mathscr{M}}
\def \N {\mathscr{N}}

\def \T {\mathscr{T}}
\def \beq {\begin{equation}}
\def \eeq {\end{equation}}
\def \ba {\begin{array}}
\def \ea {\end{array}}
\def \bs {\bigskip}

\def \ecart {\noalign{\medskip}}
\newtheorem{Thm}{Theorem}[section]
\newtheorem{Cor}[Thm]{Corollary}
\newtheorem{Pro}[Thm]{Proposition}
\newtheorem{Lem}[Thm]{Lemma}

\newtheorem{Adef}[Thm]{Definition}
\newenvironment{Def}{\begin{Adef}\rm}{\end{Adef}}
\newtheorem{Arem}[Thm]{Remark}
\newenvironment{Rem}{\begin{Arem}\rm}{\end{Arem}}
\newtheorem{Aexa}[Thm]{Example}
\newenvironment{Exa}{\begin{Aexa}\rm}{\end{Aexa}}
\newtheorem{Anot}[Thm]{Notation}

\def \refe #1.{\eqref{#1})}
\def \reff #1.{figure~\ref{#1}}
\def \refs #1.{Section~\ref{#1}}
\def \refss #1.{Subsection~\ref{#1}}
\def \refD #1.{Definition~\ref{#1}}
\def \refT #1.{Theorem~\ref{#1}}
\def \refL #1.{Lemma~\ref{#1}}
\def \refC #1.{Corollary~\ref{#1}}
\def \refP #1.{Proposition~\ref{#1}}
\def \refPt #1.{Properties~\ref{#1}}
\def \refR #1.{Remark~\ref{#1}}
\def \refE #1.{Example~\ref{#1}}
\def \refN #1.{Notation~\ref{#1}}

\newcounter{marnote}

\title{Increase of mass and nonlocal effects in the homogenization of magneto-elastodynamics problems}
\author{
\footnotesize
\centerline{\begin{tabular}{cccc}
\normalsize Marc BRIANE && \normalsize Juan CASADO-D\'IAZ
\\
Institut de Recherche Math\'ematique de Rennes && Dpto. de Ecuaciones Diferenciales y An\'alisis Num\'erico
\\
Univ Rennes, INSA de Rennes, CNRS, IRMAR - UMR 6625 && Universidad de Sevilla
\\
mbriane@insa-rennes.fr && jcasadod@us.es
\end{tabular}}
}
\begin{document}
\maketitle
\begin{abstract}
The paper deals with the homogenization of a magneto-elastodynamics equation satisfied by the displacement $u_\varepsilon$ of an elastic body which is subjected to an oscillating magnetic field $B_\varepsilon$ generating the Lorentz force $\partial_t u_\varepsilon\times B_\varepsilon$.
When the magnetic field $B_\varepsilon$ only depends on time or on space, the oscillations of $B_\varepsilon$ induce an increase of mass in the homogenized equation. More generally, when the magnetic field is time-space dependent through a uniformly bounded component $G_\varepsilon(t,x)$ of $B_\varepsilon$, besides the increase of mass the homogenized equation involves the more intricate limit $g$ of $\partial_t u_\varepsilon\times G_\varepsilon$ which turns out to be decomposed in two terms. The first term of $g$ can be regarded as a nonlocal Lorentz force the range of which is limited to a light cone at each point $(t,x)$. The cone angle is determined by the maximal velocity defined as the square root of the ratio between the elasticity tensor spectral radius and the body mass. Otherwise, the second term of $g$ is locally controlled in $L^2$-norm by the compactness default measure of the oscillating initial energy.
\end{abstract}
\section{Introduction}
In a insulating (vacuum-like) environment, an elastic three-dimensional body placed in an electric field $E$ and a magnetic $B$ is subjected to the Lorentz force (see, {\em e.g.}, \cite[Section~9.3]{BaFiFi})
\beq\label{fL}
f_L = \rho_e(E+v\times B)+\sigma(E+v\times B)\times B,
\eeq
where $v$ is the velocity, $\sigma$ the conductivity of the body and $\rho_e$ is the density of free electrical charges, while $E$ and $B$ satisfy Maxwell's system. In particular, the fields $E$ and $B$ are connected by the equation
\[
{\rm curl}\,E+\partial_t B = 0.
\]
In the present paper we focus on the magnetic Lorentz force $\partial_t v\times B$ rather than on the electrical force. We assume that
\begin{itemize}
\item the elastic body is a poor conductor, {\em i.e.} $\sigma\approx 0$,
\item the electrical Lorentz force $\rho_e E$ is negligible compared to the magnetic Lorentz force $\rho_e(v\times B)$,
\end{itemize}
which yields
\beq\label{fLs}
f_L\approx \rho_e(v\times B).
\eeq
The second assumption holds in particular if $E(t,x)=\ep\,\textsc{e}(t,x/\ep)$ with $0<\ep\ll 1$, since then
\[
O(\ep)=E(t,x) \ll B(t,x) = B(0,x)-\int_0^t ({\rm curl}\,\textsc{e})(s,x/\ep)\,ds=O(1).
\]
\par
Under these assumptions and setting $\rho_e=1$, the displacement $u$ of the body with velocity $v=\partial_t u$, satisfies the ``simplified" magneto-elastodynamics equation
\beq\label{meeq}
\rho\,\partial^2_{tt}u-{\rm Div}_x\big(A e(u)\big)+\partial_t u\times B=f,
\eeq
where $\rho$ is the mass density, $A$ is the elasticity tensor of the body and $e(u)$ is the symmetric strain tensor. The right-hand side $f$ encompasses all other body forces. Equation \eqref{meeq} can be extended to any dimension $N\geq 2$, replacing the three-dimensional Lorentz force $\partial_t u\times B$ by $B\partial_t u$, where $B$ is now a $N\times N$ skew-symmetric matrix-valued function.
\par\bs
In the framework of homogenization theory, our aim is to study the effect of a time-space oscillating magnetic field $B_\ep(t,x)$ on the magneto-elastodynamics equation \eqref{meeq}.
\par
Let $T>0$, let $\Om$ be a bounded open set of $\RR^N$ and $Q:=(0,T)\times\Om$. Consider the magneto-elastodynamics problem
\beq\label{mepb}
\left\{\ba{l}\dis \rho\,\partial^2_{tt}u_\ep-{\rm Div}_x\big(Ae(u_\ep)\big)+B_\ep\partial_tu_\ep=f_\ep\ \hbox{ in }Q
\\ \ecart
\dis u_\ep=0\ \hbox{ on }(0,T)\times\partial\Om
\\ \ecart
\dis u_\ep(0,.)=u_\ep^0,\ \partial_t u_\ep(0,.)=u_\ep^1\ \mbox{ in }\Om,
\ea\right.
\eeq
where $B_\ep$ is a skew-symmetric matrix-valued function in $L^\infty(Q)^{N\times N}$ decomposed as
\beq\label{BFGH}
B_\ep(t,x)=F_\ep(x)+G_\ep(t,x)+H_\ep(t,x),\quad\mbox{with}\quad B_\ep(0,x)=F_\ep(x),\ G_\ep(0,x)=H_\ep(0,x)=0,
\eeq
$f_\ep\in L^1(0,T;L^2(\Om))^N$, $u_\ep^0\in H^1_0(\Om)^N$, $u_\ep^1\in L^2(\Om)^N$.
Contrary to $F_\ep(x)$ the component $G_\ep(t,x)$ is assumed to be uniformly bounded with respect to $t$ and $x$, but the time-space oscillations of $G_\ep(t,x)$ may produce a nonlocal effect. The component $H_\ep(t,x)$ is a compact perturbation of $B_\ep(t,x)$.
Under suitable oscillations of the sequences $B_\ep$, $f_\ep$, $u_\ep^0$, $u_\ep^1$, we can pass to the limit as $\ep$ tends to zero in \eqref{mepb} in order to derive the homogenized problem.
\par\bs
Homogenization of pde's with variable coefficients has been the subject of numerous studies for bounded coefficients in the books \cite{Tar3,Eva,Dal,BrDe} and the references therein, \cite{BLP,San,BP} for periodically oscillating coefficients, as well as for non-uniformly bounded coefficients (from above or below) in the survey article \cite{Khr}. Here, the asymptotic analysis of equation~\eqref{mepb} involving the general sequence $B_\ep$ is in line with homogenization of pde's with non-periodic coefficients which was initiated by Spagnolo \cite{Spa} and Murat, Tartar~\cite{MuTa}. More specifically, in the stationary case Tartar \cite{Tar1,Tar2} (see also \cite{BrGe} for an alternative approach) has studied the homogenization of the three-dimensional Stokes equation
\beq\label{Seq}
-\,\De u_\ep+b_\ep\times u_\ep+\nabla p_\ep=f\ \mbox{ in }\Om,
\eeq
perturbed by the oscillating drift term $b_\ep\times u_\ep$ representing the Coriolis force which plays an analogous role to the Lorentz force \eqref{fLs} in  equation \eqref{meeq}. To that end Tartar developed his celebrated ``oscillating test functions method" at the end of the Seventies, and he obtained a homogenized Brinkman \cite{Bri} type equation
\beq\label{hSeq}
-\,\De u+b\times u+\nabla p+M^*u=f\ \mbox{ in }\Om,
\eeq
where $M^*$ is a non-negative symmetric matrix-valued function. If the magnetic field $B_\ep(x)$ is independent on time and $T=\infty$, a time Laplace transform of equation \eqref{mepb} leads us to an equation which is similar to \eqref{Seq}. Therefore, Tartar's homogenization result combined with an inverse Laplace transform should at the least modify the mass $\rho$ in the homogenized equation of~\eqref{mepb}.
\par
Alternatively, nonlocal effects without change of mass have been obtained in \cite{CCMM2} for the homogenization of a scalar wave equation with a periodically oscillating matrix-valued function $B_\ep(t,x)=B(t,x,t/\ep,x/\ep)$, where $B(t,x,s,y)$ is bounded with respect to the variables $(t,x)$ and periodically continuous with respect to the variables $(s,y)$, using a two-scale analysis method.
\par\bs
In our non-periodic and vectorial setting we show that the time-space oscillations of the magnetic field $B_\ep(t,x)$ produce both an increase of mass and nonlocal effects through an abstract representation formula arising in the homogenized equation.
\par
On the one hand, the first result of the paper is the derivation of an anisotropic effective mass $\varrho^*$ which is greater (in the sense of the quadratic forms) than the starting mass $\rho I_N$. This increase of mass in the homogenization process is due to the oscillations of the magnetic field at the microscopic scale, which modify the linear momentum through the magnetic Lorentz force. At this point Milton and Willis \cite{MiWi} have explained the macroscopic change of mass obtained in composite elastic bodies at fixed frequency, by the existence of a hidden mass at the microscopic scale, which modifies Newton's second law. From this observation, when the magnetic field $B_\ep$ is only time dependent, we can build an anisotropic internal mass $m_\ep(t)$ such that in a multiplicative way
\beq\label{mept}
m_\ep(t)\,\partial_t u_\ep\approx m^*(t)\,\partial_t u.
\eeq
In contrast, when the magnetic field is independent of time, {\em i.e.} $B_\ep=F_\ep$ which is assumed to converge weakly to zero in $W^{-1,p}(\Om)^{N\times N}$ for some $p>N$, we can build an anisotropic internal mass $M_\ep(x)=F_\ep(x) W_\ep(x)$ such that in an additive way
\beq\label{Mepx}
\partial_t u_\ep\approx\partial_t u+W_\ep\,\partial^2_{tt}u.
\eeq
The harmonic limit of $m_\ep(t)$ (due to the multiplicativity of \eqref{mept}) or the arithmetic limit of $M_\ep(x)$ (due to the additivity of \eqref{Mepx}) leads us to the anisotropic effective mass~$\varrho^*$.
\par
On the other hand, the second result of the paper shows that both time and space oscillations of the magnetic field $B_\ep(t,x)$ may also induce nonlocal effects which are absent if the magnetic field is only time dependent or only space dependent. Assuming that the component $G_\ep$ of the magnetic field (see \eqref{BFGH}) weakly converges to zero in $L^{\infty}(Q)^{N\times N}$ and that $H_\ep$ is a compact perturbation, we prove (see Theorem~\ref{thmS3}, Theorem~\ref{thmg} and Theorem~\ref{thmrep}) that the limit $g$ of the magnetic Lorentz force $G_\ep\partial_t u_\ep$ admits the following decomposition
\beq\label{decompg}
g=\int_{Q}d\Lambda(s,y)\,\partial_t u(s,y)-h_0,\ \mbox{ in  } Q.
\eeq
First, the matrix-valued measure $\Lambda$ in \eqref{decompg} can be regarded as the kernel of a nonlocal Lorentz force arising in the homogenized problem. The range of this nonlocal term is limited to each light cone of $Q$, the angle of which is equal to $2\arctan c$ with $c=\sqrt{\|A\|/\rho}$  ($\|A\|$ is the Frobenius norm of tensor $A$).
A particular nonlocal term with some light cone range has been obtained in the periodically oscillating case of \cite{CCMM2}. The general nonlocal term \eqref{decompg} with the specific light cone range is thus substantial to the homogenization process in the present non-periodic setting.
Next, we show (see Theorem~\ref{thmg}, Corollary~\ref{Thbcin} and Corollary~\ref{Thbcinf}) that the second term $h_0$ in \eqref{decompg} is locally controlled in $L^2$-norm by the compactness default measure $\mu^0$ of the oscillating initial energy. The function $h_0$ acts as a new exterior force in the homogenized problem.
\par
Therefore, collecting the two previous results we get that the homogenized problem of~\eqref{mepb} can be written as
\beq\label{hmepb}
\left\{\ba{l}
\dis \varrho^*\partial^2_{tt}u-{\rm Div}_x\big(Ae(u)\big)+H \partial_tu+g=f\ \hbox{ in }Q
\\ \ecart
\dis u=0\ \hbox{ on }(0,T)\times\partial\Om
\\ \ecart
\dis u(0,.)=u^0\ \mbox{ in }\Om,
\ea\right.
\eeq
and the initial velocity $\partial_t u(0,.)$ actually depends on the effective mass $\varrho^*$.
As a by-product of the energy estimate satisfied by the limit $g$, we obtain a corrector result for the homogenization problem \eqref{mepb}  if the compactness default measure $\mu^0$ vanishes (see Remark~\ref{remthmS4}). This holds in particular when the initial conditions are ``well-prepared" (see Remark~\ref{remarkmu0=0}) in the spirit of the classical homogenization result \cite{FrMu} for the wave equation.
\par\bs
The paper is organized as follows:
\par
In Section~\ref{S1} we study the case where the magnetic field $B_\ep$ only depends on time. We derive (see Theorem~\ref{thmS1}) the homogenized problem \eqref{hmepb} with the sole increase of mass ($g=0$).
\par
Section~\ref{S2} is devoted to a stationary problem (see Theorem~\ref{thmS2}) which prepares the main homogenization result of the paper in Section~\ref{S3}.
It is partly based on Tartar's works \cite{Tar1,Tar2}.
\par
In Section~\ref{S3} we consider a more general magnetic field $B_\ep$ satisfying \eqref{BFGH}.
We prove (see Theorem~\ref{thmS3}) that the homogenized magneto-elastodynamics problem of \eqref{mepb} is \eqref{hmepb}.
\par
Section~\ref{S4} deals with several estimates of the limit $g$ (see Theorem~\ref{thmS4} and Theorem~\ref{thmg}) and an abstract representation (see Theorem~\ref{thmrep}) which allow us to prove that the function $g$ admits the decomposition \eqref{decompg}. 
In some specific cases we get a complete representation of the function $g$ and the uniqueness of a solution to the limit problem \eqref{hmepb} (see Corollary~\ref{Thbcin} and Corollary~\ref{Thbcinf}).
\subsubsection*{Notation}
\begin{itemize}
\item $\bar Y$ denotes the closure of a subset $Y$ of a topological set $X$.
\item $A\in\L(\RR^{N\times N}_s;\RR^{N\times N}_s)$ is a positive definite symmetric fourth-order tensor, and $\|A\|$ denotes its Frobenius norm.
\item $|E|$ denotes the Lebesgue measure of a measurable set $E$ of $\RR^N$.
\item $\L(X;Y)$ denotes the space of continuous linear functions from the normed space $X$ into the normed space $Y$.
\item $\cdot$ denotes the scalar product in $\RR^N$, $:$ denotes the scalar product in $\RR^{N\times N}$, and $|\cdot|$ denotes the associated norm in both cases.
\item $B(x,r)$ denotes the euclidien ball of center $x\in\RR^N$ and of radius $r>0$, and $B_r$ simply denotes the ball $B(0,r)$ centered at the origin.
\item $I_N$ denotes the unit matrix of $\RR^{N\times N}$, and $M^t$ denotes the transposed of a matrix of $M$.
\item $\RR^{N\times N}_s$ denotes the space of symmetric matrices of order $N$.
\item $\Om$ denotes a bounded open set of $\RR^N$ for $N\geq 2$, $T>0$, and $Q$ the cylinder $(0,T)\times\Om$.
\item {\rm Div} denotes the vector-valued divergence operator taking the divergence of each row of a matrix-valued function.
\item $e(u)$ denotes the symmetrized gradient of a vector-valued function $u$.
\item $\M(X)$ denotes the space of the Radon measures on a locally compact set $X$.
\item $C^\infty_c(U)$ denotes the set of the smooth functions with compact support in an open subset $U$ of $\RR^N$.
\item $\D'(U)$ denotes the space of the distributions on an open subset $U$ of $\RR^N$.
\item $\to$ denotes a strong convergence, $\rightharpoonup$ a weak convergence, and $\stackrel{*}{\rightharpoonup}$ a weak-$*$ convergence
\item $\hookrightarrow$ denotes a continuous embedding between two topological spaces.
\item $O_\ep$ denotes a sequence of $\ep$ which converges to zero as $\ep$ tends to zero, and which may vary from line to line.
\item $C$ denotes a positive constant which may vary from line to line.
\end{itemize}
\section{Homogenization of an elastodynamics problem with a strong magnetic field only depending on time}\label{S1}
Let $\Om$ be a bounded open set of $\RR^N$ with $N\geq 2$, $T>0$,  $Q=(0,T)\times \Om$, $\rho>0$ and let $A\in\L(\RR^{N\times N}_s;\RR^{N\times N}_s)$ be a positive definite symmetric tensor. Let $f\in L^1(0,T;L^2(\Om))^N$ and let $\mbox{\ss}_\ep\in C^1([0,T])^{N\times N}$ be such that
\beq\label{Bept}
\mbox{\ss}_\ep^{t}=-\,\mbox{\ss}_\ep\;\;\mbox{and}\;\; \mbox{\ss}'_\ep \mbox{\ss}_\ep=\mbox{\ss}_\ep \mbox{\ss}'_\ep\ \mbox{ in }[0,T],
\quad  \mbox{\ss}_\ep(0)=0,
\eeq
\beq\label{cvBept}
\exp(-\rho^{-1}\mbox{\ss}_\ep)\stackrel{\ast}\rightharpoonup\mathtt{M}^{-1}\ \mbox{ in }L^\infty(0,T)^{N\times N}\ \mbox{ with }\mathtt{M}\in C^1([0,T])^{N\times N}
\eeq
where $\mbox{\ss}'_\ep$ denotes the time derivative of $\mbox{\ss}_\ep$ and $\mathtt{M}$ is an invertible matrix-valued function. Here we need to assume that the weak limit of  $\exp(-\rho^{-1}\mbox{\ss}_\ep)$ is an invertible matrix-valued function. This is not automatically satisfied as shows the following example.
\begin{Exa}\label{exa.invert}
Let $N=2$ and $\rho=1$. Consider the function $b_\ep$ defined by $b_\ep(t)=t/\ep+b(t/\ep)$ for $t\in\RR$, where $b$ is a $2\pi$-periodic function in $C^1(\RR)$, and the skew-symmetric matrix-valued function
\[
\mbox{\ss}_\ep:=\begin{pmatrix} 0 & -b_\ep \\ b_\ep & 0 \end{pmatrix}\in C^1(\RR)^{2\times 2}.
\]
We have
\[
\ba{c}
\dis \exp(-\rho^{-1}\mbox{\ss}_\ep)=\begin{pmatrix} \cos(b_\ep) & -\sin(b_\ep) \\ \sin(b_\ep) & \cos(b_\ep) \end{pmatrix}\;\rightharpoonup\;
\\ \ecart
\dis {1\over 2\pi}\int_0^{2\pi}\begin{pmatrix} \cos(t)\cos(b(t))-\sin(t)\sin(b(t)) & -\cos(t)\sin(b(t))-\sin(t)\cos(b(t)) \\ \cos(t)\sin(b(t))+\sin(t)\cos(b(t)) & \cos(t)\cos(b(t))-\sin(t)\sin(b(t))\end{pmatrix}dt.
\ea
\]
Hence, if $b=0$ then the weak limit of $ \exp(-\rho^{-1}\mbox{\ss}_\ep)$ is the nul matrix. Otherwise, if $b$ is closed to the $2\pi$-periodic function which agrees in $[0,2\pi]$ with ${\pi\over 2}\,1_{[0,\pi]}$, then the weak limit of $ \exp(-\rho^{-1}\mbox{\ss}_\ep)$ is closed to the matrix
\[
-{1\over\pi}\begin{pmatrix} 1 & -1 \\ 1 & 1 \end{pmatrix},
\]
and is thus invertible.
\end{Exa}
\par
We consider the solution $u_\ep$ to the wave equation
\beq\label{pbot}
\left\{\ba{l}
\dis \rho\,\partial^2_{tt}u_\ep-{\rm Div}_x\big(Ae(u_\ep)\big)+\mbox{\ss}'_\ep\partial_tu_\ep=f\ \hbox{ in }Q
\\ \ecart
\dis u_\ep=0\ \hbox{ on }(0,T)\times\partial\Om
\\ \ecart
\dis u_\ep(0,.)= u_\ep^0,\ \partial_t u_\ep(0,.)= u_\ep^1\ \mbox{ in }\Om,
\ea\right.
\eeq
where
\beq\label{cvcoin}
u^0_\ep\rightharpoonup u^0\ \hbox{ in }H^1_0(\Om)^N,\quad u_\ep^1\rightharpoonup u^1\ \hbox{ in }L^2(\Om)^N.
\eeq
We have the following homogenization result.
\begin{Thm}\label{thmS1}
Assume that conditions \eqref{Bept}, \eqref{cvBept}, \eqref{cvcoin} hold true. Then, we have
\beq\label{cvuep}
u_\ep\stackrel{\ast}\rightharpoonup u\ \hbox{ in }L^\infty(0,T;H^1_0(\Om))^N\cap W^{1,\infty}(0,T;L^2(\Om)))^N,
\eeq
where $u$ is the solution to the equation
\beq\label{pbothom}
\left\{\ba{l}\dis \rho\,\mathtt{M}^{t}\mathtt{M}\,\partial^2_{tt}u-{\rm Div}_x\big(Ae(u)\big)+\rho\,\mathtt{M}^{t}\mathtt{M}'\,\partial_tu=f\ \hbox{ in }Q
\\ \ecart
\dis u=0\ \hbox{ on }(0,T)\times\partial\Om
\\ \ecart
\dis u(0,.)=u^0,\ \partial_t u(0,.)=\mathtt{M}^{-1}(0)\,u^1\ \mbox{ in }\Om.
\ea\right.
\eeq
\end{Thm}
\begin{Rem}\label{rem.invert}
Since the matrix $\exp(\rho^{-1}\mbox{\ss}_\ep)$ is unitary, the lower semi-continuity of convex functionals yields for any $\la\in\RR^N$ and for any measurable set $E\subset\Om$,
\[
\int_E\kern -1.mm\mathtt{M}^{-1}\la\cdot\mathtt{M}^{-1}\la\,dx\leq\liminf_{\ep\to 0}\kern -1.mm\int_E\kern -1.mm\exp(-\rho^{-1}\mbox{\ss}_\ep)\la\cdot\exp(-\rho^{-1}\mbox{\ss}_\ep)\la\,dx=|E|\,|\la|^2,
\]
which implies that $\mathtt{M}^{-1}\la\cdot\mathtt{M}^{-1}\la\leq|\la|^2$ a.e. in $\Om$. Due to $\mathtt{M}=(\mathtt{M}^{-1})^{-1}$, we equivalently get that for any $\la\in\RR^N$, $\mathtt{M}\la\cdot\mathtt{M}\la\geq|\la|^2$ a.e. in $\Om$.
Hence, the homogenized equation~\eqref{pbothom} involves an effective anisotropic mass
\[
\rho\,\mathtt{M}^{t}\mathtt{M}\geq\rho\,I_N\ \mbox{ a.e. in }\Om,
\]
which is greater than the initial one $\rho$. We will see in Section~\ref{S3} that if we replace time oscillations by space oscillations, the homogenization process also induces a larger effective anisotropic mass but in quite a different way.
\end{Rem}
\noindent
{\bf Proof of Theorem~\ref{thmS1}.}
First of all, since $\mbox{\ss}_\ep$ is skew-symmetric, the classical estimates for the wave equation yield convergence \eqref{cvuep} up to a subsequence.
\par
By \eqref{Bept} we have $\big(\exp(\rho^{-1}\mbox{\ss}_\ep)\big)'=\rho^{-1}\,\mbox{\ss}'_\ep\exp(\rho^{-1}\mbox{\ss}_\ep)=\rho^{-1}\exp(\rho^{-1}\mbox{\ss}_\ep)\,\mbox{\ss}'_\ep$. Hence, equation~\eqref{pbot} can be written as
\beq\label{eqotvf}
\rho\,\partial_{t}\big(\exp(\rho^{-1}\mbox{\ss}_\ep)\,\partial_t u_\ep\big)-\exp(\rho^{-1}\mbox{\ss}_\ep)\,{\rm Div}_x\big(Ae(u_\ep)\big)=\exp(\rho^{-1}\mbox{\ss}_\ep)f\ \hbox{ in }Q,
\eeq
which implies that for any $\Phi\in C^1([0,T];C^\infty_c(\Om))^N$ with $\Phi(T,.)=0$,  recalling $\mbox{\ss}_\ep(0)=0$,
\beq\label{pbotvf}
\ba{ll}
\dis \int_Q\Big(\!-\!\rho\exp(\rho^{-1}\mbox{\ss}_\ep)\,\partial_t u_\ep\cdot\partial_t\Phi+Ae(u_\ep):e\big(\exp(-\rho^{-1}\mbox{\ss}_\ep)\Phi\big)\Big)\,dt\,dx
\\ \ecart
\dis =\int_\Om \rho\, u_\ep^1\cdot\Phi(0,.)\,dx+\int_Q f\cdot\,\exp(-\rho^{-1}\mbox{\ss}_\ep)\Phi\,dt\,dx.
\ea
\eeq
For $\ph\in C^\infty_c(\Om)^N$, define the function $\xi_\ep\in L^\infty(0,T)^N$ by
\beq\label{xiep}
\xi_\ep(t):=\rho\exp\big(\rho^{-1}\mbox{\ss}_\ep(t)\big)\int_\Om\partial_t u_\ep(t,x)\cdot\ph(x)\,dx,\ \mbox{ for a.e. }t\in(0,T).
\eeq
By \eqref{eqotvf} and \eqref{cvuep} we have
\[
\xi'_\ep=\int_\Om\big(\!-\!Ae(u_\ep):e(\exp(-\rho^{-1}\mbox{\ss}_\ep)\ph)+f\cdot\,\exp(-\rho^{-1}\mbox{\ss}_\ep)\ph\big)\,dx\quad\mbox{bounded in }L^\infty(0,T).
\]
Hence, $\xi_\ep$ is bounded in $W^{1,\infty}(0,T)$, and up to a subsequence converges weakly-$*$ to some $\xi$ in $W^{1,\infty}(0,T)^N$.
This combined with convergences \eqref{cvBept} and \eqref{cvuep} implies that
\[
\exp(-\rho^{-1}\mbox{\ss}_\ep)\,\xi_\ep\stackrel{*}{\rightharpoonup}\mathtt{M}^{-1}\xi=\rho\int_\Om\partial_t u(t,x)\cdot\ph(x)\,dx\ \mbox{ in }L^\infty(0,T)^N.
\]
Due to the arbitrariness of $\ph$ it follows that
\beq\label{cvBepuep}
\exp(\rho^{-1}\mbox{\ss}_\ep)\,\partial_t u_\ep\stackrel{*}{\rightharpoonup}\mathtt{M}\,\partial_t u\ \mbox{ in }L^\infty(0,T;L^2(\Om))^N.
\eeq
Moreover, integrating by parts with respect to $x$ and noting that $\mbox{\ss}_\ep$ is independent of $x$, the weak convergence \eqref{cvuep} of $u_\ep$ (see, {\em e.g.}, \cite[Chapter~3, Section~8]{LiMae}) and \eqref{cvBept} yield
\[
\ba{l}
\dis \int_Q Ae(u_\ep):e\big(\exp(-\rho^{-1}\mbox{\ss}_\ep)\Phi\big)\,dt\,dx=-\int_Q u_\ep:{\rm Div}_x\big[Ae\big(\exp(-\rho^{-1}\mbox{\ss}_\ep)\Phi\big)\big]\,dt\,dx
\\ \ecart
\dis -\int_Q u:{\rm Div}_x\big[Ae(\mathtt{M}^{-1}\Phi)\big]\,dt\,dx+O_\ep=\int_Q Ae(u):e(\mathtt{M}^{-1}\Phi)\,dt\,dx+O_\ep.
\ea
\]
Therefore, passing to the limit in \eqref{pbotvf} with \eqref{cvBepuep} and \eqref{cvuep}, we get that for any $\Phi\in C^\infty_c(Q)^N$,
\[
\int_Q\left(-\rho\,\mathtt{M}\,\partial_t u\cdot\partial_t\Phi+Ae(u):e(\mathtt{M}^{-1}\Phi)\right)dt\,dx
=\int_Q f\cdot\,\mathtt{M}^{-1}\Phi\,dt\,dx.
\]
which is equivalent to the first equation of \eqref{pbothom}.
\par
Finally, let $\Phi\in C^1([0,T];C^\infty_c(\Om))^N$ with $\Phi(T,.)=0$. Passing to the limit  in \eqref{pbotvf} with the second convergences of \eqref{cvBept} and \eqref{cvcoin} we get that
\[
\ba{ll}
\dis \int_Q\left(-\rho\,\mathtt{M}\,\partial_t u\cdot\partial_t\Phi+Ae(u):e(\mathtt{M}^{-1}\Phi)\right)dt\,dx
\\ \ecart
\dis =\int_\Om\rho\,u^1\cdot\Phi(0,.)\,dx+\int_Q f\cdot\,\mathtt{M}^{-1}\Phi\,dt\,dx,
\ea
\]
which combined with the first equation of \eqref{pbothom} gives the initial condition $\mathtt{M}(0)\,\partial_t u(0,.)=u^1$. The condition $u(0,.)=u^0$ just follows from \eqref{cvuep}, which also implies that $u_\ep$ converges to $u$ in $C^0([0,T];L^2(\Om)^N)$.
The proof is now complete.
\cqfd
\section{Homogenization of a stationary problem}
\label{S2}
Let $\Om$ be a smooth bounded open set of $\RR^N$ with $N\geq 2$.
Consider a sequence $F_\ep$ of matrix-valued functions in $W^{-1,p}(\Om)^{N\times N}$ with $p>N$, such that
\beq\label{convFep} F_\ep\rightharpoonup 0\ \hbox{ in }W^{-1,p}(\Om)^{N\times N},\eeq
and define $w_\ep^j$, $1\leq j\leq N$, as the solution to
\beq\label{ecuacwei} \left\{\ba{l}\dis -\,{\rm Div}\big(Ae(w_\ep^j)\big)+F_\ep e_j=0\ \hbox{ in }\Om\\ \ecart\dis w_\ep^j=0\ \hbox{ on }\partial\Om,\ea\right.\eeq
which satisfies (due to the regularity of $\Om$)
\beq\label{conwei}
w_\ep^j\rightharpoonup 0\ \hbox{ in }W^{1,p}_0(\Om)^N,\quad\forall\,j\in \{1\ldots,N\}.
\eeq
Extracting a subsequence if necessary, we can assume the existence of a nonnegative symmetric matrix-valued function $M$ in $ L^{p\over 2}(\Om)^{N\times N}$ such that
\beq\label{defM} 
Ae(w_\ep^j):e(w_\ep^k)\rightharpoonup (Me_j)\cdot e_k\ \mbox{ in }L^{p\over 2}(\Om),\quad\forall\,j,k\in\{1,\dots,N\}.
\eeq
\par
We have the following result which will be used in the next section with $u_\ep$ as a time average of the displacement and $z_\ep$ as a time average of the velocity in the elastodynamics problem.
\begin{Thm} \label{thmS2}
Consider two sequences $z_\ep\in H^1(\Om)^N$ and $f_\ep\in H^{-1}(\Om)^N$ such that
\beq\label{conzefe} z_\ep\rightharpoonup z\ \hbox{ in }H^1(\Om)^N,\qquad f_\ep\rightarrow f\ \hbox{ in }H^{-1}(\Om)^N,\eeq
and define $u_\ep$ as the solution to
\beq\label{ecuuep}\left\{\ba{l}\dis -\,{\rm Div}\big(Ae(u_\ep)\big)+F_\ep z_\ep=f_\ep\ \hbox{ in }\Om\\ \ecart\dis u_\ep=0\ \hbox{ on }\partial\Om.\ea\right.\eeq
Then, up to a subsequence, we have
\beq\label{conuep} u_\ep\rightharpoonup u\ \hbox{ in }H^1_0(\Om)^N\eeq
\beq\label{corresc}u_\ep-u-\sum_{j=1}^Nw^j_\ep z_j\to 0\ \hbox{ in }H^1_0(\Om)^N,\eeq
\beq\label{objet}
Ae(u_\ep):e(u_\ep)\stackrel{\ast}\rightharpoonup Ae(u):e(u)+Mz\cdot z\ \hbox{ in }\M(\bar\Om).
\eeq
\end{Thm}\noindent
{\bf Proof.}
First of all, observe that thanks to Rellich-Kondrachov's compactness theorem, we have
\beq\label{DEDEDE}
v_\ep^1,v_\ep^2\rightharpoonup 0\ \hbox{ in }H^1(\Om)^N\;\;\Rightarrow\;\; F_\ep v_\ep^1\cdot v_\ep^2\rightharpoonup 0\ \hbox{ in }\D'(\Om).
\eeq
Using that $F_\ep={\rm div}\,(\Phi_\ep)$ where $\Phi_\ep$ is bounded in $L^p(\Om)^N$, and that by Sobolev's imbedding $\nabla(z_\ep u_\ep)$ belongs to $L^{N\over N-1}(\Om)^N\subset L^{p'}(\Om)^N$ (due to $p>N$) where $z_\ep$ is bounded in $H^1(\Om)$, we have by H\"older's inequality
\[
\big|\langle F_\ep,z_\ep u_\ep\rangle\big|\leq \|\Phi_\ep\|_{L^p(\Om)^N}\,\|\nabla (z_\ep u_\ep)\|_{L^{N\over N-1}(\Om)^N}
\leq C\,\|u_\ep\|_{H^1_0(\Om)}.
\]
Hence, putting $u_\ep$ as test function in \eqref{ecuuep} the former estimate combined with the boundedness of $f_\ep$ in $H^{-1}(\Om)$
implies that $u_\ep$ is bounded in $H^1_0(\Om)$. Therefore, convergence \eqref{conuep} holds up to a subsequence.
\par
Now,
given $\phi\in C^\infty_c( \Om)^N$, we put 
$$u_\ep-u-\sum_{j=1}^Nw^j_\ep\phi_j$$ 
as test function in \eqref{ecuuep}. Thanks to \eqref{DEDEDE} and \eqref{conwei}, 
we get
\beq\label{primec}\into Ae(u_\ep):\Big(e(u_\ep)-e(u)-\sum_{j=1}^Ne(w^j_\ep)\phi_j\Big)\,dx +\left\langle F_\ep z,\Big(u_\ep-u-\sum_{j=1}^Nw^j_\ep\phi_j\Big)\right\rangle
=O_\ep.
\eeq
On the other hand, putting 
\[
\phi_j\Big(u_\ep-u-\sum_{i=1}^Nw^i_\ep\phi_i\Big)
\]
as test function in \eqref{ecuacwei}, adding in $j$ and using \eqref{conwei}, we get
\beq\label{secec}\ba{l}\dis \into A\Big(\sum_{j=1}^Ne(w_\ep^j)\phi_j\Big):\Big(e(u_\ep)-e(u)-\sum_{j=1}^Ne(w^j_\ep)\phi_j\Big)\,dx\\ \ecart\dis+\left\langle F_\ep\phi,(u_\ep-u-\sum_{j=1}^Nw^j_\ep\phi_j\Big )\right\rangle=O_\ep.
\ea\eeq
Subtracting \eqref{primec} and  \eqref{secec} we have
\[
\ba{l}
\dis\into A\Big(e(u_\ep)-\sum_{j=1}^Ne(w^j_\ep)\phi_j\Big ):\Big(e(u_\ep)-e(u)-\sum_{j=1}^Ne(w^j_\ep)\phi_j\Big )\,dx
\\ \ecart
\dis +\left\langle F_\ep(z-\phi),\Big(u_\ep-u-\sum_{j=1}^Nw^j_\ep\phi_j\Big)\right\rangle=O_\ep.
\ea
\]
This combined with the weak convergence
\[
e(u_\ep)-e(u)-\sum_{j=1}^Ne(w^j_\ep)\phi_j\;\rightharpoonup 0\quad\mbox{in }L^2(\Om)^{N\times N}
\]
also yields
\[
\ba{l}
\dis\into A\Big(e(u_\ep)-e(u)-\sum_{j=1}^Ne(w^j_\ep)\phi_j\Big ):\Big(e(u_\ep)-e(u)-\sum_{j=1}^Ne(w^j_\ep)\phi_j\Big )\,dx
\\ \ecart
\dis +\left\langle F_\ep(z-\phi),\Big(u_\ep-u-\sum_{j=1}^Nw^j_\ep\phi_j\Big)\right\rangle=O_\ep.
\ea
\]
From Rellich-Kondrachov's compactness theorem and Cauchy-Schwarz' inequality, we deduce
\[
\limsup_{\ep\to 0}\into A\Big(e(u_\ep)-e(u)-\sum_{j=1}^Ne(w^j_\ep)\phi_j\Big ):\Big(e(u_\ep)-e(u)-\sum_{j=1}^Ne(w^j_\ep)\phi_j\Big )\,dx
\leq C\,\|z-\phi\|_{H^1_0(\Om)^N}.
\]
Moreover, taking a sequence $\phi^n$ which converges strongly to $z$ in $H^1_0(\Om)^N$ and noting that 
\[
\lim_{n\to\infty}\limsup_{\ep\to 0}\into A\Big(\sum_{j=1}^Ne(w_\ep^j)(\phi_j^n-z_j)\Big):\Big(\sum_{j=1}^Ne(w_\ep^j)(\phi_j^n-z_j)\Big)\,dx=0,
\]
we conclude to 
\beq\label{corres}
\lim_{\ep\to 0}\into A\Big(e(u_\ep)-e(u)-\sum_{j=1}^Ne(w^j_\ep)z_j\Big):\Big(e(u_\ep)-e(u)-\sum_{j=1}^Ne(w^j_\ep)z_j\Big)\,dx=0,
\eeq
which by Korn's inequality proves \eqref{corresc}.
\par\noindent
It is immediate that \eqref{corres} and \eqref{defM} imply \eqref{objet}.\cqfd
\par\medskip
We also have the following lower semicontinuity result.
\par
\begin{Lem} \label{lemscps}
Consider a sequence $u_\ep$ which satisfies the assumptions of Theorem~\ref{thmS2}.
Then, up to subsequence, there exists a measurable function $\zeta:\Om\to\RR^N$, with
\beq\label{proprlFu}
M\zeta\cdot \zeta\in L^1(\Om)^N,\quad M\zeta\in L^{2p\over p+2}(\Om)^N,
\eeq
such that
\beq\label{conFeFue}
F^{t}_\ep u_\ep\rightharpoonup M \zeta\ \hbox{ in }H^{-1}(\Om)^N,\eeq
\beq\label{semicFeu}
\liminf_{n\to\infty}\into Ae(u_\ep):e(u_\ep)\varphi\,dx\geq \into \big(Ae(u):e(u)+M\zeta\cdot \zeta\big)\varphi\,dx,\quad \forall\, \varphi\in C^0(\bar\Om),\ \varphi\geq 0.
\eeq
\end{Lem}
\noindent
{\bf Proof.} For $\phi\in C^\infty(\bar \Om)^N$, we have
\[
\big\langle F_\ep^{t} u_\ep,\phi\big\rangle=\big\langle F_\ep\phi,u_\ep\big\rangle=-\sum_{j=1}^N\into Ae(u_\ep ):e(w_\ep^j)\phi_j\,dx+O_\ep.
\]
Therefore, defining $Z\in L^{2p\over p+2}(\Om)^N$ by
\beq\label{elesc0}
Ae(u_\ep):e(w_\ep^j)\rightharpoonup -Z_j\ \hbox{ in }L^{2p\over p+2}(\Om),
\eeq
we get
\beq\label{elesc1}
F_\ep^{t} u_\ep\rightharpoonup Z\ \hbox{ in }H^{-1}(\Om)^N.
\eeq
\par
Applying \eqref{elesc1} to $w_\ep^k$ in place of $u_\ep$ and recalling the definition of $M$, we get 
\beq\label{elesc2}
F_\ep^{t}w_\ep^k\rightharpoonup -Me_k\ \hbox{ in }H^{-1}(\Om)^N,
\eeq
which implies
\beq\label{elesc3}
F_\ep^{t}\Big(\sum_{k=1}^Nw_\ep^k\phi_k\Big)\rightharpoonup -M\phi\ \hbox{ in }H^{-1}(\Om)^N,\quad \forall\, \phi\in C^\infty(\bar\Om)^N.
\eeq
\par
On the other hand, by \eqref{elesc0}, \eqref{defM} and Cauchy-Schwarz' inequality, we have for any function $\eta\in L^{2p\over p-2}(\Om)^N$,
\[
\ba{l}
\dis \left|\,\into Z\cdot\eta\,dx\,\right|=\left|\,\lim_{\ep\to 0}\into Ae(u_\ep):\Big(\sum_{j=1}^Ne(w_\ep^j)\eta_j\Big)\,dx\,\right|
\\ \ecart
\dis \leq \left(\limsup_{\ep\to 0}\into Ae(u_\ep):e(u_\ep)\,dx\right)^{1\over 2}\left(\into M\eta\cdot\eta\,dx\right)^{1\over 2}.
\ea
\]
Therefore, $Z$ is orthogonal to any function $\eta\in L^{2p\over p-2}(\Om)^N$ such that $M\eta=0$ a.e. in $\Om.$
Now, let us show the existence of a measurable function $\zeta:\Om\to\RR^N$ such that $Z=M\zeta.$
To this end, consider the set
\[
V:=\left\{M\xi:\xi\in L^{2p\over p-1}(\Om)^N\right\}
\]
which by H\"older's inequality (recall that $M\in L^{p\over 2}(\Om)^{N\times N}$) is a linear subspace of $L^{2p\over p+1}(\Om)^N$.
Let $\eta\in L^{2p\over p-2}(\Om)^N$. Due to the symmetry of $M$ we have
\[
\eta\in V^\perp\;\Leftrightarrow\;\forall\,\xi\in L^{2p\over p+1}(\Om)^N,\ \int_\Om M\eta\cdot\xi\,dx=0\;\Leftrightarrow\; M\eta=0\mbox{ a.e. in }\Om,
\]
which implies that $Z\in(V^\perp)^\perp=\overline{V}$ since $L^{2p\over p+1}(\Om)^N$ is a reflexive space (see, {\em e.g.}, \cite[Proposition~1.9]{Bre}). Hence, there exists a sequence $\zeta_n$ in $L^{2p\over p-1}(\Om)^N$ such that $M\zeta_n$ converges strongly to $Z$ in $L^{2p\over p+1}(\Om)^N$. Up to replace $\zeta_n$ by its orthogonal projection on $(\ker M)^\perp$, we may assume that $\zeta_n\in(\ker M)^\perp$ a.e. in $\Om$.
Next, consider the measurable pseudo-inverse $M^{-1}$ of the matrix-valued $M$ defined for a.e. $x\in\Om$ by $M(x)^{-1}(M(x)\,\xi)=\xi$ for any $\xi\in (\ker M(x))^\perp$. Then, we have for any $k>0$,
\[
1_{\{|M^{-1}|\leq k\}}\,\zeta_n=1_{\{|M^{-1}|\leq k\}}M^{-1}(M\zeta_n)\to 1_{\{|M^{-1}|\leq k\}}M^{-1}Z\quad\mbox{strongly in }L^{2p\over p+1}(\Om)^N.
\]
Since a strongly convergent sequence in $L^q(\Om)$ converges up to a subsequence a.e. in $\Om$, using a diagonal procedure in the former convergences there exists a subsequence $\zeta_{\theta(n)}$ such that for any $k>0$,
\[
1_{\{|M^{-1}|\leq k\}}\,\zeta_{\theta(n)}\to 1_{\{|M^{-1}|\leq k\}}\,M^{-1}Z\quad\mbox{a.e. in }\Om.
\]
Therefore, the a.e. limit $\zeta$ of $\zeta_{\theta(n)}$ in $\Om$ satisfies $Z=M\zeta$ a.e. in $\Om$.
\par
It remains to prove \eqref{semicFeu} which in particular implies the first assertion of \eqref{proprlFu}.
Taking into account \eqref{conwei} and \eqref{conuep}, for any $\phi\in C^0( \bar \Om)^N$ and $\varphi\in C^0(\bar\Om)$, $\varphi\geq 0$, we have
\[
\ba{l}
\dis \into A\Big(e(u_\ep)-e(u)-\sum_{j=1}^Ne(w^j_\ep)\,\phi_j\Big ):\Big(e(u_\ep)-e(u)-\sum_{j=1}^Ne(w^j_\ep)\,\phi_j\Big)\,\varphi\,dx
\\ \ecart
\dis =\into Ae(u_\ep):e(u_\ep)\,\varphi\,dx-\into Ae(u):e(u)\,\varphi\,dx
\\ \ecart
\dis +\into A\Big(\sum_{j=1}^Ne(w^j_\ep)\,\phi_j\Big):\Big(\sum_{j=1}^Ne(w^j_\ep)\,\phi_j\Big)\varphi\,dx
 -2\sum_{j=1}^N\into Ae(u_\ep):e(w_\ep^j)\,\phi_j\,\varphi\,dx+O_\ep
\\ \ecart
\dis =\into Ae(u_\ep):e(u_\ep)\,\varphi\,dx-\into Ae(u):e(u)\,\varphi\,dx+\into M\phi:\phi\,\varphi\,dx+2\into M\zeta\cdot \phi\,\varphi\,dx.
\ea
\]
This proves 
\[
\lim_{\ep\to 0}\into Ae(u_\ep):e(u_\ep)\,\varphi\,dx\geq \into Ae(u):e(u)\varphi\,dx-\into M\phi\cdot \phi\,\varphi\,dx-2\into M\zeta\cdot\phi\,\varphi\,dx,
\]
for any $\phi\in C^0(\bar \Om)^N$ and any $\varphi\in C^0(\bar\Om)$, $\varphi\geq 0$. Taking into account that $M\zeta$ belongs to $L^{2p\over p+2}(\Om)^N$, we deduce by approximation that the above equality holds for any $\phi\in L^{2p\over p-2}(\Om)^N$. Thus we can choose in particular
$\phi=-\,1_{B(0,R)\cap \{|\zeta|< R\}}\,\zeta$.
Then, passing to the limit as $R$ tends to infinity thanks to the monotone convergence theorem we conclude to \eqref{semicFeu}. As a by-product we deduce from  \eqref{semicFeu} with $\ph=1$ that $M\zeta\cdot\zeta\in L^1(\Om)$.
\cqfd
\section{Homogenization of a general magneto-elastodynamics problem}\label{S3}
Let $\Om$ be a smooth bounded open set of $\RR^N$, $N\geq 2$, $T>0$,  $Q=(0,T)\times \Om$, $\rho>0$ and let $A\in\L(\RR^{N\times N}_s;\RR^{N\times N}_s)$ be a positive definite symmetric tensor.
\par
Consider a sequence $F_\ep$ of skew-symmetric matrix-valued functions in $L^\infty(\Om)^{N\times N}$  which satisfies \eqref{convFep} for some $p>N$, a sequence of skew-symmetric matrix-valued functions $G_\ep$ in $L^\infty(Q)^{N\times N}$ such that
\beq\label{conGep}
G_\ep\stackrel{\ast}\rightharpoonup 0\ \hbox{ in }L^\infty(Q)^{N\times N},
\eeq
and a sequence $H_\ep$ of skew-symmetric matrix-valued functions in $L^\infty(Q)^{N\times N}$ such that 
\beq\label{convHep}
H_\ep\rightarrow H\ \hbox{ in }H^1(0,T;W^{-1,p}(\Om))^{N\times N}.
\eeq
Define 
\beq\label{deBep}
B_\ep(t,x):=F_\ep(x)+G_\ep(t,x)+H_\ep(t,x)\quad \mbox{for }(t,x)\in Q.
\eeq
Recall that $M$ is the non-negative symmetric matrix-valued function in $L^{p\over 2}(\Om)^{N\times N}$ defined by~\eqref{defM}.
\par
The main result of the section is the following
\begin{Thm}\label{thmS3}
Let $f_\ep\in L^1(0,T;L^2(\Om))^N$ be such that 
\beq\label{confep} f_\ep\rightarrow f\ \hbox{ in }L^1(0,T;L^2(\Om))^N,\eeq
and $u^0_\ep\in H^1_0(\Om)^N$, $u^1_\ep\in L^2(\Om)^N$ be such that
\beq\label{concoin}
u^0_\ep\rightharpoonup u^0\ \hbox{ in }H^1_0(\Om)^N,\quad
u_\ep^1\rightharpoonup u^1\ \hbox{ in }L^2(\Om)^N.
\eeq
Then, there exist a measurable function $\zeta:\Om\to\RR^N$ and a function $g\in L^\infty(0,T;L^2(\Om))^N$ such that the solution $u_\ep$ of
\beq\label{pb1}
\left\{\ba{l}\dis \rho\,\partial^2_{tt}u_\ep-{\rm Div}_x\big(Ae(u_\ep)\big)+B_\ep \partial_tu_\ep=f_\ep\ \hbox{ in }Q
\\ \ecart
\dis u_\ep=0\ \hbox{ on }(0,T)\times\partial\Om
\\ \ecart
\dis u_\ep(0,.)=u_\ep^0,\ \partial_t u_\ep(0,.)=u_\ep^1\ \mbox{ in }\Om,
\ea\right.
\eeq
and $u^0_\ep$ satisfy up to a subsequence
\beq\label{convuep}
u_\ep\stackrel{\ast}\rightharpoonup u\ \hbox{ in }L^\infty(0,T;H^1_0(\Om))^N\cap W^{1,\infty}(0,T;L^2(\Om))^N,
\eeq
\beq\label{Defig}
G_\ep\partial_tu_\ep\stackrel{\ast}\rightharpoonup g\ \hbox{ in }L^\infty(0,T;L^2(\Om))^N,
\eeq
\beq\label{Defzet}
F_\ep u_\ep^0\rightharpoonup M\zeta \hbox{ in }H^{-1}(\Om)^N,\ \hbox{ with } M\zeta\in L^{2p\over p+2}(\Om)^N,\ M\zeta\cdot\zeta\in L^1(\Om).
\eeq
Moreover, the limit $u$ is a solution to
\beq\label{pb2}
\left\{\ba{l}
\dis (\rho I_N+M)\partial^2_{tt}u-{\rm Div}_x\big(Ae(u)\big)+H \partial_tu+g=f\ \hbox{ in }Q
\\ \ecart
\dis u=0\ \hbox{ on }(0,T)\times\partial\Om
\\ \ecart
\dis u(0,.)=u^0,\ \partial_t u(0,.)=(\rho I_N+M)^{-1}(\rho u^1+M\zeta)\ \mbox{ in }\Om,
\ea\right.
\eeq
with
\beq\label{regdtu}
M\partial_tu\cdot\partial_t u\in L^\infty(0,T;L^1(\Om)).
\eeq
\end{Thm}
\begin{Rem}\label{rem.thmS3}
Actually, the function $g$ given by convergence \eqref{Defig} is independent of the sequence $H_\ep$ (which is in some sense compact) but cannot be determined in terms of the limits $f$, $u^0$, $u^1$.
In particular, we cannot prove an uniqueness result for the limit problem~\eqref{pb2}.
In Section~\ref{S4} we will give a specific representation about the function $g$ illuminating possible nonlocal effects in the homogenization process.
\par
However, if $\partial_t G_\ep$ is assumed for instance to be bounded in $L^1(0,T;L^\infty(\Om))^{N\times N}$ which corresponds to the absence of time oscillations, then the function $g$ is zero. In this case the limit problem \eqref{pb2} is completely determined and has a unique solution. The limit elastodynamics equation~\eqref{pb2} is then characterized by a magnetic field $H$ and an increase of mass $M$ which only depends on the space oscillations of $F_\ep(x)$ through \eqref{defM}. This completes the picture of Section~\ref{S1} where the magnetic field only depends on time. The general case with both space and time oscillations through $G_\ep(t,x)$ is much more intricate and leads to the undetermined function $g$.
\par
Note that the strong convergence \eqref{convHep} makes $H_\ep$ a compact perturbation of the magnetic field which simply gives the limit $H$ in the homogenized equation \eqref{pb2}.
\end{Rem}
\noindent
{\bf Proof of Theorem~\ref{thmS3}.}  First of all (see, {\em e.g.} \cite[Chapter 1]{LiMa}), it is classical that the limit problem \eqref{pb1} has one solution in $C^0([0,T]; H^1_0(\Om))^N\cap C^1([0,T];L^2(\Om))^N$ and that, taking into account that $F_\ep$, $H_\ep$ are skew-symmetric, we have the energy identity
\[
{1\over 2}{d\over dt}\left(\into \Big(\rho|\partial_t u_\ep|^2+Ae(u_\ep):e(u_\ep)\Big)\,dx\right)=\into f_\ep\cdot\partial_t u_\ep\,dx.
\]
This implies that up to a subsequence $u_\ep$ satisfies \eqref{convuep}, \eqref{Defig}. In particular, we have
\beq\label{conpt1}
u_\ep\rightarrow u\ \hbox{ in }C^0([0,T]; L^2(\Om))^N.
\eeq
Moreover, we recall that $u\in L^\infty(0,T; H^1_0(\Om))^N\cap W^{1,\infty}(0,T;L^2(\Om))^N$ gives 
\[
\left\{\ba{l}\dis u(t,.)\in H^1_0(\Om)^N,\quad\forall\, t\in [0,T]\\ \ecart\dis
t_n\to t\;\;\Rightarrow\;\; u(t_n,.)\rightharpoonup u(t,.)\ \hbox{ in }H^1_0(\Om)^N,\ea\right.
\]
and that \eqref{convuep} implies
\beq\label{contunt}
u_\ep(t,.)\rightharpoonup u(t,.)\ \hbox{ in }H^1_0(\Om)^N,\quad \forall\, t\in [0,T].\eeq
\par
Now, the idea is to take time-average values of $u_\ep$ and to apply the results of Section \ref{S2}.
Integrating  \eqref{pb1} with respect to $t$ in $(t_1,t_2)$ with $0\leq  t_1<t_2\leq T$, we deduce that the function
\[
\bar u_\ep:=\int_{t_1}^{t_2} u_\ep(s,.)\,ds\ \hbox{ in }\Om
\]
satisfies
\beq\label{ecuaint}\ba{l}\dis \rho\big(\partial_tu_\ep(t_2,x)-\partial_tu_\ep(t_1,x)\big)-{\rm Div}_x \big(Ae(\bar u_\ep)\big)+F_\ep\big(u_\ep(t_2,x)-u_\ep(t_1,x)\big)
\\ \ecart
\dis +(H_\ep u_\ep)(t_2,x)-(H_\ep u_\ep)(t_1,x)-\int_{t_1}^{t_2}\partial_tH_\ep u_\ep\,dt+\int_{t_1}^{t_2}G_\ep\partial_t u_\ep\,dt
\\ \ecart
\dis=\int_{t_1}^{t_2}f_\ep\,dt\ \hbox{ in }H^{-1}(\Om)^N.\ea
\eeq
\noindent
{\it First step.} A corrector result for $\bar u_\ep$.
\par\noindent
By \eqref{convuep} we have 
\beq\label{conpt}
\rho\big(\partial_tu_\ep(t_2,.)-\partial_tu_\ep(t_1,.)\big)\ \hbox{ bounded in }L^2(\Om)^N.
\eeq
By \eqref{convHep} we have
\[
H_\ep(t,.)\rightarrow H(t,.)\ \hbox{ in }C^0([0,T];W^{-1,p}(\Om))^{N\times N},
\]
which combined with  \eqref{contunt} and Rellich-Kondrachov's theorem gives
\beq\label{constt}(H_\ep u_\ep)(t_2,.)-(H_\ep u_\ep)(t_1,.)\rightarrow (H u)(t_2,.)-(H u)(t_1,.)\ \hbox{ in }H^{-1}(\Om)^N.\eeq
Similarly, we have
\beq\label{conct}\int_{t_1}^{t_2}\partial_tH_\ep u_\ep\,dt\to \int_{t_1}^{t_2}\partial_tH u\,dt\ \hbox{ in }H^{-1}(\Om)^N.\eeq
By \eqref{Defig} we also have
\beq\label{consst}\int_{t_1}^{t_2}G_\ep\partial_t u_\ep\,dt\rightharpoonup \int_{t_1}^{t_2}g\,dt\ \hbox{ in }L^2(\Om)^N.\eeq
The previous convergences \eqref{conpt}, \eqref{constt}, \eqref{conct}, and \eqref{consst} combined with  \eqref{contunt} and \eqref{ecuaint} allow us to apply Theorem \ref{thmS2} to deduce the corrector result
\beq\label{corresCDT}
\bar u_\ep-\bar u-\sum_{j=1}^N\big(u_j(t_2,.)-u_j(t_1,.)\big)w_\ep^j\to 0\ \hbox{ in }H^1(\Om)^N,
\eeq 
where we denote
\[
\bar u:=\int_{t_1}^{t_2} u(s,.)\,ds\ \hbox{ in }\Om.
\]
\par\medskip\noindent
{\it Second step.} Limit of \eqref{ecuaint}.
\par\noindent
We replace in \eqref{ecuaint}, $t_1$, $t_2$ by $t_1+s$, $t_2+s$ and we integrate with respect to $s$ in $(0,\tau)$ with $\tau<T-t_2$. Then, we can pass to the limit as $\ep$ tends to zero to deduce
\beq\label{ecuaintt}
\ba{l}
\dis \rho\big(u(t_2+\tau,.)-u(t_2,.)-u(t_1+\tau,.)+u(t_1,.)\big)
-{\rm Div} \left(Ae\Big(\int_0^\tau\int_{t_1+s}^{t_2+s}u\,dt\,ds\Big)\right)
\\ \ecart
\dis+\lim_{\ep\to 0}F_\ep\left(\int_0^\tau\big(u_\ep(t_2+s,.)-u_\ep(t_1+s,.)\big)\,ds\right) +\int_0^\tau\big((H u)(t_2+s,.)-(H u)(t_1+s,.)\big)\,ds
\\ \ecart
\dis +\int_0^\tau\int_{t_1+s}^{t_2+s}\big(-\partial_tH u+g\big)\,dt\,ds=\int_0^\tau\int_{t_1+s}^{t_2+s}f\,dt\,ds\ \hbox{ in }H^{-1}(\Om)^N,\ea\eeq
where the limit in the third term is taken in the weak topology of $H^{-1}(\Om)^N$. Moreover, by \eqref{convFep}, \eqref{corresCDT},  \eqref{elesc2} and $F_\ep$  skew-symmetric we have
\beq\label{limFeuep}
\ba{l}
\dis F_\ep\left(\int_0^\tau\big(u_\ep(t_2+s,.)-u_\ep(t_1+s,.)\big)\,ds\right)
=\dis F_\ep\left(\int_{t_2}^{t_2+\tau}u_\ep(s,.)\,ds-\int_{t_1}^{t_1+\tau}u_\ep(s,.)\,ds\right)
\\ \ecart
\dis =\sum_{j=1}^N \big(u_j(t_2+\tau,.)-u_j(t_2,.)-u_j(t_1+\tau,.)+u_j(t_1,.)\big)F_\ep w^j_\ep+R_\ep
\\ \ecart
\dis =M\big(u(t_2+\tau,.)-u(t_2,.)-u(t_1+\tau,.)+u(t_1,.)\big)+R_\ep,
\ea
\eeq
where $R_\ep$ denotes a sequence which converges weakly (strongly for the first one) to zero in $H^{-1}(\Om)^N$.
Putting  \eqref{limFeuep} in \eqref{ecuaintt}, dividing by $\tau$ and letting $\tau$ tend to zero, we get
$$\ba{l}\dis \rho\big(\partial_tu(t_2,.)-\partial_tu(t_1,.)\big)-{\rm Div}\left(Ae\Big(\int_{t_1}^{t_2}u\,dt\Big)\right)\\ \ecart\dis+M\big(\partial_tu(t_2,.)-\partial_tu(t_1,.)\big)
+(H u)(t_2,.)-(H u)(t_1,.)\\ \ecart\dis+
\int_{t_1}^{t_2}\big(-\partial_tH u+g\big)\,dt=\int_{t_1}^{t_2}f\,dt\ \hbox{ in }H^{-1}(\Om)^N.\ea$$
Finally, dividing by $t_2-t_1$ and letting $t_2-t_1$ tend to zero, we obtain
\beq\label{Primecli}(\rho I_N+M)\partial^2_{tt}u-{\rm Div}_x\big(Ae(u)\big)+H \partial_tu+g=f\ \hbox{ in }\D'(Q)^N.\eeq
\par\medskip\noindent
{\it Third step.} Limit of the initial conditions.
\par\noindent
By \eqref{contunt} and \eqref{concoin} the limit $u$ satisfies
\beq\label{condinu}
u(0,.)=u^0\ \hbox{ in }\Om.
\eeq
\par
Now, it remains to find the initial velocity. Let us prove that
\beq\label{condere}
(\rho\partial_tu_\ep+F_\ep u_\ep)(t,.)\rightharpoonup \big((\rho I_N+M)\partial_tu\big)(t,.)\ \hbox{ in }H^{-1}(\Om)^N,\quad\forall\, t\in [0,T).
\eeq
By \eqref{pb1} we have
\[
\partial_t\big(\rho\partial_tu_\ep+F_\ep u_\ep+H_\ep u_\ep\big)=f_\ep+{\rm Div}_x\big(Ae(u_\ep)\big)+\partial_tH_\ep u_\ep-G_\ep\partial_tu_\ep,
\]
where the right-hand side is bounded in $L^1(0,T;H^{-1}(\Om))^N$ by \eqref{conGep}, \eqref{convHep}, \eqref{confep}, \eqref{convuep}.
Therefore, 
$$\rho\partial_tu_\ep+F_\ep u_\ep+H_\ep u_\ep\ \hbox{ is bounded in }W^{1,1}(0,T;H^{-1}(\Om))^N,$$
Now, we fix $t_0\in [0,T)$ and we observe that \eqref{convFep}, $u_\ep\in C^0([0,T];H^1_0(\Om))^N\cap C^1([0,T];L^2(\Omega))^N$ and \eqref{convuep} imply, up to a subsequence,
\beq\label{conciL} (\rho\partial_tu_\ep+F_\ep u_\ep)(t_0,.)\rightharpoonup L\ \hbox{ in } H^{-1}(\Om)^N.\eeq
On the other hand, for $\tau\in (0,T-t_0)$, we have
\[
\ba{l}\dis\left\|\,(\rho\partial_tu_\ep+F_\ep u_\ep+H_\ep u_\ep)(t_0,.)-{1\over \tau}\int_{t_0}^{t_0+\tau}(\rho\partial_tu_\ep+F_\ep u_\ep+H_\ep u_\ep)(t,.)\,dt\,\right\|_{H^{-1}(\Om)^N}
\\ \ecart
\dis =\left\|\,{1\over \tau}\int_0^\tau\left(\int_{t_0}^{t_0+t}  \partial_r(\rho\partial_r u_\ep+F_\ep u_\ep+H_\ep u_\ep)(r,.)\,dr\right)\,dt\,\right\|_{H^{-1}(\Om)^N}
\\ \ecart
\dis \leq {1\over \tau}\int_0^\tau\left(\int_{t_0}^{t_0+t} \big\|f_\ep+{\rm Div}_x\big(Ae(u_\ep)\big)+\partial_r H_\ep u_\ep-G_\ep\partial_r u_\ep\big\|_{H^{-1}(\Om)^N}\,dr\right)\,dt
\\ \ecart
\dis \leq \|f_\ep(t,.)\|_{L^1(t_0,t_0+\tau;L^2(\Om))^N}+\left({\tau\over 2}+
C\sqrt{\tau}\,\big\|\partial_tH_\ep \big\|_{L^2(t_0,t_0+\tau;W^{-1,p}(\Om)^{N\times N})}\right)\|u_\ep\|_{L^\infty(0,T;H^1_0(\Om))^N}
\\ \ecart
\dis +\,C\tau\,\|G_\ep\|_{L^\infty(Q)^N}\|\partial_tu_\ep\|_{L^\infty(0,T;L^2(\Om))^N}.
\ea
\]
By \eqref{convHep}, \eqref{convuep} and \eqref{contunt} we have
\[
(H_\ep u_\ep)(t_0,.)\rightarrow (H u)(t_0,.)\ \hbox{ in }H^{-1}(\Om)^N,
\]
\[
{1\over \tau}\int_{t_0}^{t_0+\tau}(\rho\partial_tu_\ep+H_\ep u_\ep)(t,.)\,dt\rightarrow
{1\over \tau}\int_{t_0}^{t_0+\tau}(\rho\partial_tu+H u)(t,.)\,dt\ \hbox{ in }H^{-1}(\Om)^N.
\]
By \eqref{limFeuep} we also have
\[
{1\over \tau}\int_{t_0}^{t_0+\tau}F_\ep u_\ep\,dt\rightharpoonup {1\over \tau}\int_{t_0}^{t_0+\tau} M\partial_t u\,dt\ \hbox{ in }H^{-1}(\Om).
\]
Therefore, we deduce
\beq\label{ecu1Lecoi}
\ba{l}\dis\left\|\, L+(H u)(t_0,.)-{1\over \tau}\int_{t_0}^{t_0+\tau}\big((\rho I_N+M)\partial_tu+Hu\big)(t,.)\,dt\,\right\|_{H^{-1}(\Om)^N}
\\ \ecart
\dis \leq \|f(t,.)\|_{L^1(t_0,t_0+\tau;L^2(\Om))^N}+C\big(\sqrt{\tau}+\tau\big).
\ea
\eeq
Next, equation \eqref{Primecli}, combined with $u\in L^\infty(0,T;H^1_0(\Om))^N\cap W^{1,\infty}(0,T;L^2(\Om))^N$ implies as above
\[
(\rho I_N+M)\partial_tu+Hu\in W^{1,1}(0,T; H^{-1}(\Om))^N\hookrightarrow C^0([0,T];H^{-1}(\Om))^N.
\]
Hence, passing to the limit in \eqref{ecu1Lecoi} as $\tau$ tends to zero, we get
\[
L=\big((\rho I_N+M)\partial_tu\big)(t_0,.),
\]
which implies \eqref{condere}.
\par\medskip
Convergence \eqref{condere} combined with \eqref{concoin} yields
\beq\label{convconInd}
\rho u_\ep^1+F_\ep u_\ep^0\rightharpoonup \big((\rho I_N+M)\partial_tu\big)(0,.)\ \hbox{ in }H^{-1}(\Om)^N.
\eeq
Therefore, by \eqref{conFeFue} and $F_\ep$ skew-symmetric there exists a measurable function $\zeta$ satisfying \eqref{Defzet}, which yields the second initial condition of \eqref{pb2}.
\par
Finally, the proof of estimate \eqref{regdtu} is given in Lemma~\ref{lemsecef} below. This concludes the proof of Theorem~\ref{thmS3}.
\cqfd
\section{Energy estimates and  nonlocal effects}\label{S4}
The aim of this section is to estimate more precisely the function $g$ arising in the homogenized problem~\eqref{pb2}.
\subsection{Energy estimate}
\par\noindent
First of all, observe that the following inequality holds
\beq\label{ineM}
(\rho I_N+M)^{-1}(\rho\xi+M\eta)\cdot(\rho\xi+M\eta)\leq \rho|\xi|^2+M\eta\cdot\eta,\quad\forall\,\xi,\eta\in\RR^N.
\eeq
In order to show it, set $\Upsilon:=(\rho I_N+M)^{-1}(\rho\xi+M\eta)$. Then, using successively the Cauchy-Schwarz inequality with the non-negative symmetric matrix $M$ and the Cauchy-Schwarz inequality in $\RR^2$, we have 
\[
(\rho I_N+M)\Upsilon\cdot \Upsilon=(\rho\xi+M\eta)\cdot\Upsilon\leq (\rho|\xi|^2+M\eta\cdot\eta)^{1\over 2}\big((\rho I_N+M)\Upsilon\cdot \Upsilon\big)^{1\over 2},
\]
which gives \eqref{ineM}.
\par
From {\eqref{concoin}, applying the lower semicontinuity \eqref{semicFeu} and convergence~\eqref{conFeFue} with $u_\ep=u_\ep^0$, and applying the inequality \eqref{ineM} with $\xi=u^1$ and $\eta=\zeta$, we can assume, up to extract a subsequence, that there exists a non-negative Radon measure $\mu^0$ defined on $\bar \Om$ such that
\beq\label{wepococ}
\ba{ll}
& \dis \rho|u_\ep^1|^2+ Ae(u_\ep^0):e(u_\ep^0)
\\ \ecart
\stackrel{\ast}\rightharpoonup & \dis \mu^0+Ae(u^0):e(u^0)+(\rho I_N+M)^{-1}(\rho u^1+M\zeta)\cdot(\rho u^1+M\zeta)\ \hbox{ in }\M(\bar\Om).
\ea
\eeq
\begin{Rem}\label{remarkmu0=0}
The measure $\mu^0$ represents the compactness default with respect to the initial conditions $u_\ep^0$, $u_\ep^1$. Now, assume that the initial conditions are well-prepared (see \cite{FrMu} for the classical homogenization of the wave equation without Lorentz force) in the following sense:
\beq\label{condbp}
u_\ep^1\rightharpoonup u^1\ \hbox{ in }H^1(\Om)^N,\quad -\,{\rm div}\big(Ae(u_\ep^0)\big)+F_\ep u_\ep^1\ \hbox{ is compact in }H^{-1}(\Om)^N.
\eeq
Then, using the convergence \eqref{corresc} with $u_\ep=u_\ep^0$ and $z_\ep=u_\ep^1$, combined with convergences \eqref{convFep}, \eqref{elesc2}, we get that $M\zeta=Mu^1$. Moreover, by \eqref{objet} and Rellich-Kondrachov's compactness theorem we have
\[
\rho|u_\ep^1|^2+ Ae(u_\ep^0):e(u_\ep^0)
\;\stackrel{\ast}\rightharpoonup\;\rho|u^1|^2+ Ae(u^0):e(u^0)+Mu^1\cdot u^1\ \hbox{ in }\M(\bar\Om),
\]
which proves  that the measure $\mu^0$ vanishes.
\end{Rem}
\par\medskip
Let us introduce the following notations.
\begin{Def}
Set
\beq\label{defvitma}
c:=\sqrt{\|A\|\over \rho}.
\eeq
For $\bar x\in \bar\Om$, $S\in (0,T)$ and $t\in (0,S)$, we denote
\beq\label{defBSs}
B(\bar x,S,t):= B\big(\bar x,c(S-t)\big)\cap \Om.
\eeq
\beq\label{dePBSs}
K(\bar x,S,t):=\partial B\big(\bar x,c(S-t)\big)\cap \Om,
\eeq
and recall that $B_\delta$ is the ball centered at the origin of radius $\delta>0$.
\par
\end{Def}
We have the following result.
\begin{Thm}\label{thmS4}
Under the assumptions and the notations of Theorem~\ref{thmS3}, for any $\bar x\in \bar\Om$,  $0<S_1<S_2$, $s\in (0,S_1)$, $\delta>0$, and $\psi\in L^2(0,T;L^{2p\over p-2}(\Om))^N$, the solution $u_\ep$ of \eqref{pb2} satisfies
\beq\label{estifu}
\ba{l}
\dis \limsup_{\ep\to 0}{1\over 2}\int_{S_1}^{S_2}\!\!\int_{B_\delta}\int_{B(\bar x+z,S,s)}\Big[\rho|\partial_t(u_\ep-u)|^2 +
\\ \ecart
\dis \kern 1cm +A\Big(e(u_\ep)-e(u)-\sum_{j=1}^Ne(w^j_\ep)\psi^j\Big):\Big(e(u_\ep)-e(u)-\sum_{j=1}^Ne(w^j_\ep)\psi^j\Big)\Big]dx\,dz\,dS
\\ \ecart
\dis \leq  {1\over 2}\int_{S_1}^{S_2}\!\!\int_{B_\delta} \mu^0\big(\bar B(\bar x+z,S,0)\big)\,dz\,dS
+\int_{S_1}^{S_2}\!\!\int_{B_\delta}\int_{B(\bar x+z,S,t)}\kern -.4em M(\partial_t u-\psi)\cdot(\partial_t u-\psi)\,dx\,dz\,dS
\\ \ecart
\dis +\int_{S_1}^{S_2}\!\!\int_{B_\delta}\int_0^s\int_{B(\bar x+z,S,t)}\kern -.4em g\cdot\partial_t u\,dx\,dt\,dz\,dS.
\ea
\eeq
\par
\end{Thm}
\begin{Rem} \label{remthmS4}
Assuming that the measure $\mu^0$ and the function $g$ vanish, estimate \eqref{estifu} gives the corrector result
\beq\label{estifu2}
\ba{l}
\dis \limsup_{\ep\to 0}{1\over 2}\int_{S_1}^{S_2}\!\!\int_{B_\delta}\int_{B(\bar x+z,S,s)}\Big[\rho|\partial_t(u_\ep-u)|^2+
\\ \ecart
\dis \kern 1cm +A\Big(e(u_\ep)-e(u)-\sum_{j=1}^Ne(w^j_\ep)\psi^j\Big):\Big(e(u_\ep)-e(u)-\sum_{j=1}^Ne(w^j_\ep)\psi^j\Big)\Big]dx\,dz\,dS
\\ \ecart
\dis \leq 
\int_{S_1}^{S_2}\!\!\int_{B_\delta}\int_{B(\bar x+z,S,t)}\kern -.4em M(\partial_t u-\psi)\cdot(\partial_t u-\psi)\,dx\,dz\,dS.
\ea
\eeq
By virtue of Remark \ref{rem.thmS3} and Remark \ref{remarkmu0=0} a sufficient condition for estimate \eqref{estifu2} to be satisfied is that the sequence $\partial_t G_\ep$ is bounded in $L^1(0,T;L^\infty(\Om))^{N\times N}$ and that the initial conditions are well-prepared in the sense of \eqref{condbp}.
\end{Rem}
From Theorem \ref{thmS4} we deduce the following estimate for the function $g$ which will be improved in Section \ref{sectiong}.
\begin{Cor}\label{estimateg}
Under the same assumptions of Theorem \ref{thmS4}, there exists a constant $C>0$ which only depends on $\sup_{\ep>0}\|G_\ep\|_{L^\infty(Q)}^{N\times N}$ such that the function $g$ of \eqref{pb2} satisfies
\beq\label{Estifg2}
\int_{B(\bar x,S,s)}|g|^2dx\leq C\mu^0\big(\bar B(\bar x,S,0)\big)+C\left(\int_0^s\left(\int_{B(\bar x,S,t)}|\partial_tu|^2dx\right)^{1\over 2}dt\right)^2,
\eeq
\beq\label{Estifg2b}
0\leq{1\over 2}\,\mu^0\big(\bar B(\bar x,S,0)\big)+\int_0^s\int_{B(\bar x,S,t)} g\cdot\partial_t u\,dx\,dt,
\eeq
for any $\bar x\in \bar\Om$, any $S\in (0,T)$ and a.e. $s\in (0,S)$.
\end{Cor}
\noindent
\begin{Rem}\label{remcone}
For $\bar x\in \bar\Om$ and $S\in (0,T)$, define the cone of vertex $(\bar x,S)$ and angle equal to $2\arctan c$
\beq\label{cone}
\C(S,\bar x):=\big\{(t,x):\ 0<t<S,\ x\in B\big(\bar x,c(S-t)\big)\big\},
\eeq
where  $c$ is the wave propagation velocity defined by \eqref{defvitma}. Then, estimate \eqref{Estifg2} means that the norm of $g$ over the cone section at time $t=s$ is bounded by the measure $\mu^0$ of the cone section at time $t=0$ plus the norm of the velocity $\partial_t u$ over the truncated cone in the time interval~$(0,s)$.
\end{Rem}
\noindent
{\bf Proof of Corollary \ref{estimateg}.}
By \eqref{Defig} and \eqref{estifu}  there exists a constant $C>0$ which only depends on $\sup\|G_\ep\|_{L^\infty(Q)}^{N\times N}$ such that for any $S_1,S_2$ with $0<S_1<S_2<T$, $s\in (0,S_1)$, $\de>0$, and $\psi\in L^2(0,T;L^{2p\over p-2}(\Om))^N$,
\beq\label{estimDtu2}
\ba{l}
\dis\int_{S_1}^{S_2}\!\!\int_{B_\delta}\int_{B(\bar x+z,S,s)}|g|^2\,dx\,dz\,dS
\leq  C\int_{S_1}^{S_2}\!\!\int_{B_\delta} \mu^0\big(\bar B(\bar x+z,S,0)\big)\,dz\,dS
\\ \ecart
\dis +\, C\int_{S_1}^{S_2}\!\!\int_{B_\delta}\int_{B(\bar x+z,S,t)}\kern -.4em M(\partial_t u-\psi)\cdot(\partial_t u-\psi)\,dx\,dz\,dS
\\ \ecart
\dis +\, C\int_{S_1}^{S_2}\!\!\int_{B_\delta}\int_0^s\int_{B(\bar x+z,S,t)}\kern -.4em g\cdot\partial_t u\,dx\,dt\,dz\,dS.
\ea
\eeq
Moreover, by virtue of \eqref{regdtu} and using an approximation by truncation in the space of the functions $v\in L^\infty(0,T;L^2(\Om))^N$ with $Mv\cdot v\in L^\infty(0,T;L^1(\Om))$, the sequence $\psi_n:=\partial_t u\,1_{\{|\partial_t u|\leq n\}}$ in $L^\infty(Q)^N$ satisfies
\[
\lim_{n\to\infty}\big\|M(\partial_t u-\psi_n)\cdot(\partial_t u-\psi_n)\big\|_{L^1(Q)}=0.
\]
Using this approximation in \eqref{estimDtu2} it follows that
\beq\label{estgmu0Dtu}
\ba{l}
\dis\int_{S_1}^{S_2}\!\!\int_{B_\delta}\int_{B(\bar x+z,S,s)}|g|^2\,dx\,dz\,dS
\leq  C\int_{S_1}^{S_2}\!\!\int_{B_\delta} \mu^0\big(\bar B(\bar x+z,S,0)\big)\,dz\,dS
\\ \ecart
\dis +\,C\int_{S_1}^{S_2}\!\!\int_{B_\delta}\int_0^s\int_{B(\bar x+z,S,t)}\kern -.4em g\cdot\partial_t u\,dx\,dt\,dz\,dS.
\ea
\eeq
Making $S_1$, $S_2$ tend to $S$, then $\delta$ tend to zero, this implies that for any $S\in (0,T)$ and a.e. $s\in (0,S)$,
\beq\label{Estifg}
\int_{B(\bar x,S,s)}|g|^2dx\leq C\mu^0\big( \bar B(\bar x,S,0)\big)+C\int_0^s\int_{B(\bar x,S,t)}\hskip-4pt  |g\cdot\partial_t u|\,dx\,dt.
\eeq
Now, defining
\[
\Phi(s):=\int_{B(\bar x,S,s)}|g|^2dx, \quad A:=C\mu^0\big(  \bar B(\bar x,S,0)\big),\quad K(s):=\left(\int_{B(\bar x,S,s)}|\partial_tu|^2dx\right)^{1\over 2},
\]
and using the Cauchy-Schwarz inequality in \eqref{Estifg}, it follows that
\[
\Phi(s)\leq A+C\int_0^sK(t)\,\Phi(t)^{1\over 2}\,dt.
\]
By a Gronwall's type argument this provides \eqref{Estifg2} for another constant $C$, which concludes the proof of \eqref{Estifg2}.
\par
The proof of \eqref{Estifg2b} easily follows from  \eqref{estgmu0Dtu} by taking $S_1=S$, dividing by $(S_2-S)\,\delta^N$, then letting this quantity tend to zero.
\cqfd
\par\medskip
To prove Theorem \ref{thmS4} we need the following results.
\begin{Lem}\label{LemDeriv}
Let $\Om$ be a smooth ($C^1$) open set in $\RR^N$ and $u\in W^{1,1}(0,T;L^1(\Om))\cap L^1(0,T;W^{1,1}(\Om))$. For $x_0\in \RR^N$ and  $R\in C^1(0,T)$, $R>0$, we define
\[
\Phi(t):=\int_{B(x_0,R(t))\cap\Om}\hskip-1.em u(t,x)\,dx,\ \mbox{ for }\,t\in [0,T].
\]
Then, $\Phi\in W^{1,1}(0,T)$ and 
\beq\label{derivinv}
\Phi'(t)=\int_{B(x_0,R(t))\cap\Om}\hskip-1.em \partial_tu(t,x)\,dx
+R'(t)\int_{\partial B(x_0,R(t))\cap\Om}\hskip-1.em u(t,x)\,ds(x),\ \hbox{ for a.e. }t\in (0,T).
\eeq
\end{Lem}
\begin{Lem} \label{lemsecef} The limit $u$ of the solution $u_\ep$ of  \eqref{pb1} satisfies \eqref{regdtu}.
\par\noindent
Moreover, for any $\nu\in C^0(\bar Q)^N$ with $|\nu|\leq 1$, and for any $\varphi\in C^0(\bar Q)$ with $\varphi\geq 0$, we have
\beq\label{semlf}
\ba{l}
\dis \liminf_{\ep\to 0}\int_0^{T}\int_\Om\Big({c\over 2}|\partial_tu_\ep|^2+{c\over 2}Ae(u_\ep):e(u_\ep)-Ae(u_\ep):\big(\partial_t u_\ep\odot \nu\big)\Big)\varphi\,dx\,dt
\\ \ecart\dis
\geq \int_0^{T}\int_\Om\Big({c\over 2}\big(\rho I_N+M)\partial_tu\cdot\partial_t u+{c\over 2}Ae(u):e(u)-Ae(u):\big(\partial_t u\odot \nu\big)\Big)\varphi\,dx\,dt.
\ea
\eeq
We also have
\beq\label{liuews}
Ae(u_\ep):e(w_\ep^j)\stackrel{\ast}\rightharpoonup M\partial_t u\cdot e_j\ \hbox{ in }L^\infty\big(0,T;L^{{2p\over p+2}}(\Om)\big),\quad\forall\, j\in \{1,\ldots,N\}.
\eeq
\end{Lem}
\noindent
{\bf Proof of Theorem \ref{thmS4}.} Let $\bar x\in\bar \Om$, $S\in (0,T)$, $t\in (0,S)$ and $\delta>0$.
\par
Note that if $\partial_t u_\ep\in L^\infty(0,T;H^1_0(\Om))^N$, we can put it as test function in \eqref{pb1}.
Then, integrate with respect to $x$ over $B(\bar x+z,S,t)$ and to $z$ over $B_\delta$, we get by virtue of Lemma~\ref{LemDeriv}
\beq\label{enereq}
\ba{l}\dis {1\over 2}{d\over dt}\int_{B_\delta}\int_{B(\bar x+z,S,t)}\kern -.4em\big(\rho|\partial_tu_\ep|^2+Ae(u_\ep):e(u_\ep)\big)\,dx\,dz
\\ \ecart
\dis=\int_{B_\delta}\int_{B(\bar x+z,S,t)}\kern -.4em\big(\rho\,\partial^2_{tt}u_\ep\cdot\partial_tu_\ep+Ae(u_\ep):e(\partial_tu_\ep)\big)\,dx\,dz
\\ \ecart
\dis
-{c\over 2}\int_{B_\delta}\int_{K(\bar x+z,S,t)}\hskip-4pt \big(\rho|\partial_tu_\ep|^2+Ae(u_\ep):e(u_\ep)\big)\,ds(x)\,dz
\\ \ecart
\dis =\int_{B_\delta}\int_{B(\bar x+z,S,t)}\kern -.4em\big(\rho\,\partial^2_{tt}u_\ep-{\rm Div}_x(Ae(u_\ep))\big)\cdot\partial_tu_\ep\,dx\,dz
\\ \ecart
\dis +\int_{B_\delta}\int_{K(\bar x+z,S,t)}\hskip-4pt\Big(Ae(u_\ep):\partial_tu_\ep\odot\nu-{c\over 2}\rho|\partial_tu_\ep|^2-{c\over 2}Ae(u_\ep):e(u_\ep)\Big)ds(x)\,dz\\ \ecart\dis =\int_{B_\delta}\int_{B(\bar x+z,S,t)}\kern -.4em f_\ep\cdot\partial_t u_\ep\,dx\,dz
\\ \ecart
\dis +\int_{B_\delta}\int_{K(\bar x+z,S,t)}\hskip-4pt\Big(Ae(u_\ep):\partial_tu_\ep\odot\nu-{c\over 2}\rho|\partial_tu_\ep|^2-{c\over 2}Ae(u_\ep):e(u_\ep)\Big)ds(x)\,dz.
\ea
\eeq
In the general case \eqref{enereq} remains true using an approximation argument.
\par
Integrating with respect to $t$ in $(0,s)$ with $0<s<S$, we obtain
\beq\label{estenint}
\ba{l}
\dis{1\over 2}\int_{B_\delta}\int_{B(\bar x+z,S,s)}\hskip-4pt\Big(\rho|\partial_tu_\ep|^2+Ae(u_\ep):e(u_\ep)\Big)\,dx\,dz
\\ \ecart
\dis -{1\over 2}\int_{B_\delta}\int_{B(\bar x+z,S,0)}\hskip-4pt \Big(\rho|u_\ep^1|^2+Ae(u^0_\ep):e(u^0_\ep)\Big)\,dx\,dz=
\int_{B_\delta}\int_0^s\int_{B(\bar x+z,S,t)}\kern -.4em f_\ep\cdot\partial_t u_\ep\,dx\,dt\,dz
\\ \ecart
\dis-\int_{B_\delta}\int_0^s\int_{K(\bar x+z,S,t)}\hskip-4pt  \Big({c\over 2}\rho|\partial_tu_\ep|^2+{c\over 2}Ae(u_\ep):e(u_\ep)-Ae(u_\ep):\partial_tu_\ep\odot\nu\Big)ds(x)\,dt\,dz.
\ea
\eeq
Using estimate \eqref{semlf} in \eqref{estenint} and recalling \eqref{wepococ} we then deduce
\beq\label{estimacot1}\ba{l}
\dis\limsup_{\ep\to 0}\left({1\over 2}\int_{B_\delta}\int_{B(\bar x+z,S,s)}\hskip-4pt \big(\rho|\partial_tu_\ep|^2+Ae(u_\ep):e(u_\ep)\big)\,dx\,dz\right)
\\ \ecart
\dis \leq \int_{B_\delta}\int_0^s\int_{B(\bar x+z,S,t)}\kern -.4em f\cdot\partial_t u\,dx\,dt\,dz
+{1\over 2}\int_{B_\delta}\mu^0\big(\bar B(\bar x+z,S,0)\big)\,dz
\\ \ecart
\dis +\,{1\over 2}\int_{B(\bar x+z,S,0)}\kern -.4em\big(Ae(u^0):e(u^0)+(\rho I_N+M)^{-1}(\rho u^1+M\zeta):(\rho u^1+M\zeta)\big)\,dx\,dz
\\ \ecart
\dis-\int_{B_\delta}\int_0^s\int_{K(\bar x+z,S,t)}\hskip-4pt \Big({c\over 2}(\rho I_N+M)\partial_tu\cdot\partial_t u+{c\over 2}Ae(u):e(u)-Ae(u):\partial_tu\odot\nu\Big)ds(x)\,dt\,dz.
\ea
\eeq
Moreover, the non-negativity of the last integral of \eqref{estenint}, convergences \eqref{confep}, \eqref{concoin}, \eqref{convuep}, and the inclusion $B(\bar x+z,S,s)\subset\Om$ imply that there exists a constant $C_\delta$ such that
\beq\label{dominener}
\int_{B_\delta}\int_{B(\bar x+z,S,s)}\hskip-4pt\Big(\rho|\partial_tu_\ep|^2+Ae(u_\ep):e(u_\ep)\Big)\,dx\,dz\leq C_\delta.
\eeq
\par
Next, similarly to \eqref{estenint} with equation \eqref{pb2} we have
\beq\label{estimacot2}
\ba{l}
\dis{1\over 2}\int_{B_\delta}\int_{ B(\bar x+z,S,s)}\hskip-4pt \big((\rho I_N+M)\partial_tu\cdot\partial_t u+Ae(u):e(u)\big)\,dx\,dz
\\ \ecart
\dis= {1\over 2}\int_{B_\delta}\int_{ B(\bar x+z,S,0)}\hskip-4pt \big(Ae(u^0):e(u^0)+(\rho I_N+M)^{-1}(\rho u^1+M\zeta)\cdot(\rho u^1+M\zeta)\big)\,dx\,dz
\\ \ecart
\dis+ \int_{B_\delta}\int_0^s\int_{B(\bar x+z,S,t)}\kern -.4em (f-g)\cdot\partial_t u\,dx\,dt\,dz
\\ \ecart
\dis-\int_{B_\delta}\int_0^s\int_{K(\bar x+z,S,t)}\hskip-4pt \Big({c\over 2}(\rho I_N+M)\partial_tu\cdot\partial_tu+{c\over 2}Ae(u):e(u)-Ae(u):\partial_tu\odot\nu\Big)ds(x)\,dt\,dz.
\ea
\eeq\par
On the other hand, for any $\psi\in C^\infty_c(Q)^N$, we have
\[
\ba{l}
\dis {1\over 2}\int_{S_1}^{S_2}\!\!\int_{B_\delta}\int_{B(\bar x+z,S,s)}\Big[\rho|\partial_t(u_\ep-u)|^2  +
\\ \ecart
\dis \kern 1cm +A\Big(e(u_\ep)-e(u)-\sum_{j=1}^Ne(w^j_\ep)\psi^j\Big):\Big(e(u_\ep)-e(u)-\sum_{j=1}^Ne(w^j_\ep)\psi^j\Big)\Big]dx\,dz\,dS
\\ \ecart
\dis = {1\over 2}\int_{S_1}^{S_2}\!\!\int_{B_\delta}  \int_{B(\bar x+z,S,s)}\hskip-4pt \left(\rho|\partial_tu_\ep|^2+ Ae(u_\ep):e(u_\ep)\right)\,dx\,dz\,dS
\\ \ecart
\dis -\int_{S_1}^{S_2}\!\!\int_{B_\delta}  \int_{B(\bar x+z,S,s)}\hskip-4pt \left(\rho\partial_tu_\ep\cdot\partial_t u+ Ae(u_\ep):\Big(e(u)+\sum_{j=1}^Ne(w^j_\ep)\psi^j\Big)\right)\,dx\,dz\,dS
\\ \ecart
\dis
+{1\over 2}\int_{S_1}^{S_2}\!\!\int_{B_\delta}  \int_{B(\bar x+z,S,s)}\hskip-4pt \left(\rho|\partial_tu|^2+ Ae(u):e(u)\right)\,dx\,dz\,dS
\\ \ecart
\dis+\int_{S_1}^{S_2}\!\!\int_{B_\delta}  \int_{B(\bar x+z,S,s)}\hskip-4pt Ae(u):\Big(\sum_{j=1}^Ne(w^j_\ep)\psi^j\Big)\,dx\,dz\,dS
\\ \ecart
\dis+{1\over 2}\int_{S_1}^{S_2}\!\!\int_{B_\delta}  \int_{B(\bar x+z,S,s)}
\hskip-4pt  A\Big(\sum_{j=1}^Ne(w^j_\ep)\psi^j\Big):\Big(\sum_{j=1}^Ne(w^j_\ep)\psi^j\Big)\,dx\,dz\,dS.
\ea
\]
Passing to the limit as $\ep$ tends to zero thanks to \eqref{liuews} and \eqref{defM} we get
\[
\ba{l}
\dis  \limsup_{\ep\to 0}{1\over 2}\int_{S_1}^{S_2}\!\!\int_{B_\delta}\int_{B(\bar x+z,S,s)}\Big[\rho|\partial_t(u_\ep-u)|^2 +
\\ \ecart
\dis \kern 1cm +\,A\Big(e(u_\ep)-e(u)-\sum_{j=1}^Ne(w^j_\ep)\psi^j\Big):\Big(e(u_\ep)-e(u)-\sum_{j=1}^Ne(w^j_\ep)\psi^j\Big)\Big]dx\,dz\,dS
\\ \ecart
\leq \dis \limsup_{\ep\to 0}\left({1\over 2}\int_{S_1}^{S_2}\!\!\int_{B_\delta}  \int_{B(\bar x+z,S,s)}\hskip-4pt \left(\rho|\partial_tu_\ep|^2+ Ae(u_\ep):e(u_\ep)\right)\,dx\,dz\,dS\right)
\\ \ecart
\dis -\,{1\over 2}\int_{S_1}^{S_2}\!\!\int_{B_\delta}  \int_{B(\bar x+z,S,s)}\hskip-4pt \left(\rho|\partial_t u|^2+ Ae(u):e(u)+ 2M\partial_t u\cdot\psi-M\psi\cdot\psi\right)\,dx\,dz\,dS,
\ea
\]
 which, by the Lebesgue dominated convergence theorem together with estimate \eqref{dominener}, yields
\beq\label{estimDtu}
\ba{l}
\dis  \limsup_{\ep\to 0}{1\over 2}\int_{S_1}^{S_2}\!\!\int_{B_\delta}\int_{B(\bar x+z,S,s)}\Big[\rho|\partial_t(u_\ep-u)|^2+
\\ \ecart
\dis \kern 1cm +\,A\Big(e(u_\ep)-e(u)-\sum_{j=1}^Ne(w^j_\ep)\psi^j\Big):\Big(e(u_\ep)-e(u)-\sum_{j=1}^Ne(w^j_\ep)\psi^j\Big)\Big]dx\,dz\,dS
\\ \ecart
\leq \dis \int_{S_1}^{S_2}\limsup_{\ep\to 0}\left({1\over 2}\int_{B_\delta}  \int_{B(\bar x+z,S,s)}\hskip-4pt \left(\rho|\partial_tu_\ep|^2+ Ae(u_\ep):e(u_\ep)\right)\,dx\,dz\right)dS
\\ \ecart
\dis -\,{1\over 2}\int_{S_1}^{S_2}\!\!\int_{B_\delta}  \int_{B(\bar x+z,S,s)}\hskip-4pt \left(\rho|\partial_t u|^2+ Ae(u):e(u)+ 2M\partial_t u\cdot\psi-M\psi\cdot\psi\right)\,dx\,dz\,dS.
\ea
\eeq
Estimate \eqref{estimDtu} combined with \eqref{estimacot1} and \eqref{estimacot2} finally yields \eqref{estifu} for $\psi\in C^\infty_c(Q)^N$.
The case  where $\psi\in L^2(0,T;L^{2p\over p-2}(\Om))^N$ easily follows by approximating $\psi$ by a sequence in $C^\infty_c(Q)^N$.
\cqfd
\subsection{Fine estimate of the function $g$}\label{sectiong}
Corollary \ref{estimateg} can be improved by the following result.
\begin{Thm}\label{thmg}
Under the assumptions of Theorem \ref{thmS3}
there exist a subsequence of $\ep$ still denoted by $\ep$, a constant $C>0$ which only depends on $\sup_{\ep>0}\|G_\ep\|_{L^\infty(Q)^{N\times N}}$ and a continuous linear operator $\G:L^1(0,T;L^2(\Om))^N\to L^\infty(0,T;L^2(\Om))^N$ such that for any $w\in L^1(0,T;L^2(\Om))^N$, any $\bar x\in \bar\Om$, any $S\in (0,T)$ and a.e. $s\in (0,S)$,
\beq\label{EstimateG1}
\int_{B(\bar x,S,s)}\big|\G w\big|^2dx\leq C\left(\int_0^s\left(\int_{B(\bar x,S,t)}|w|^2dx\right)^{1\over 2}dt\right)^2,
\eeq
\beq\label{EstimateG1b}
0\leq\int_0^s\int_{B(\bar x,S,t)} \big(\G w\big)\cdot w\,dx\,dt
\eeq
and such that the functions $g$ and $u$ in the limit problem \eqref{pb2} defined up to a subsequence of~$\ep$, satisfy
\beq\label{EstimateG2}
\int_{B(\bar x,S,s)}\big|g-\G(\partial_t u)\big|^2dx\leq C\,\mu^0\big( \bar B(\bar x,S,0)\big),
\eeq
where $\mu_0$ is the measure defined by \eqref{wepococ} up to a subsequence of $\ep$.
\end{Thm}
\begin{Rem}\label{RemThmG}
Theorem \ref{thmg} shows that the function $g$ of problem \eqref{pb2} is the difference of $\G(\partial_t u)$ and a function $h_0$ which only depends on the initial conditions $u_\ep^0$, $u_\ep^1$ of problem \eqref{pb1} through the measure $\mu^0$.  The additional term $h_0$ acts as a new exterior force in the limit equation \eqref{pb2}.
\end{Rem}
As a consequence of Theorem \ref{EstimateG1b} we can now get a full representation of the limit  problem \eqref{pb2} for some particular choices of the initial conditions. Our first result refers to the case of well-prepared initial conditions in the sense of Remark \ref{remarkmu0=0}.
\begin{Cor}\label{Thbcin}
Consider the subsequence of $\ep$ defined by Theorem \ref{thmg}. Assume that the initial conditions $u^0_\ep$, $u^1_\ep$ in \eqref{pb1} satisfy \eqref{condbp}. Then, the solution $u_\ep$ of \eqref{pb1} satisfies \eqref{convuep}, where $u$ is the unique solution to
\beq\label{pb2Cibp}
\left\{\ba{l}
\dis (\rho I_N+M)\partial^2_{tt}u-{\rm Div}_x\big(Ae(u)\big)+H \partial_tu+\G(\partial_t u)=f\ \hbox{ in }Q
\\ \ecart
\dis u=0\ \hbox{ on }(0,T)\times\partial\Om
\\ \ecart
\dis u(0,.)=u^0,\ \partial_t u(0,.)=u^1\ \mbox{ in }\Om,
\ea\right.
\eeq
\end{Cor}
 As an example of not well-prepared initial data consider the case where the initial conditions do not depend on $\ep$.
\begin{Cor}\label{Thbcinf}
There exists a subsequence of $\ep$ such that Theorem \ref{thmg} holds and such that there exists a constant $C>0$, which only depends on $\sup_{\ep>0}\|G_\ep\|_{L^\infty(Q)^N}$ and a continuous linear operator $\F:L^2(\Om)^N\to L^\infty(0,T;L^2(\Om))^N$ such that for every $v\in L^2(\Om)^N$, any $\bar x\in \bar\Om$, any $S\in (0,T)$ and a.e. $s\in (0,S)$,
\beq\label{estimopP}
\int_{ B(\bar x,S,s)}|\F (v)|^2dx\leq C\int_{ B(\bar x,S,0)}(\rho I+M)^{-1}Mv\cdot v\,dx,
\eeq
and such for any $u^0\in H^1_0(\Om)$ and $u^1\in L^2(\Om)$, the solution $u_\ep$ of \eqref{pb1} with $u^0_\ep=u^0$, $u^1_\ep=u^1$ satisfies \eqref{convuep}, where $u$ is the unique solution to 
\beq\label{pb2CipP}
\left\{\ba{l}
\dis (\rho I_N+M)\partial^2_{tt}u-{\rm Div}_x\big(Ae(u)\big)+H \partial_tu+\G(\partial_t u)=f+\F(u^1)\ \hbox{ in }Q
\\ \ecart
\dis u=0\ \hbox{ on }(0,T)\times\partial\Om
\\ \ecart
\dis u(0,.)=u^0,\ \partial_t u(0,.)=\rho (\rho I_N+M)^{-1}u^1\ \mbox{ in }\Om,
\ea\right.
\eeq
Moreover, for any $\bar x\in \bar\Om$, any $S\in (0,T)$ and a.e. $s\in (0,S)$,
\beq\label{positp}
\int_0^s\int_{B(\bar x,S,s)}\hskip-10pt\F (u^1)\cdot \partial_tu\,dx\,dt\leq \int_0^s\int_{B(\bar x,S,s)}\hskip-10pt\G(\partial_tu)\cdot\partial_tu\,dx\,dt+{\rho\over 2}\int_{ B(\bar x,S,0)}\hskip-4pt(\rho I+M)^{-1}Mu^1\cdot u^1\,dx.
\eeq
\end{Cor}
The proof of Theorem \ref{thmg} is based on the following result.
\begin{Lem} \label{lemthmg}
Let $w\in C^\infty_c(Q)^N$ and let $v_\ep^k$ for $k\in\NN$ be the solution to
\beq\label{pbaw}
\left\{\ba{l}\dis \rho\,\partial^2_{tt}v^k_\ep-{\rm Div}_x\big(Ae(v^k_\ep)\big)+(F_\ep+G_\ep) \partial_tv^k_\ep+k\big(\partial_t v^k_\ep-w\big)=0\ \hbox{ in }Q
\\ \ecart
\dis v^k_\ep=0\ \hbox{ on }(0,T)\times\partial\Om
\\ \ecart
\dis v^k_\ep(0,.)=0,\ \partial_t v^k_\ep(0,.)=0\ \mbox{ in }\Om.
\ea\right.
\eeq
Then, there exists a constant $C_w>0$ such that for any $k\in \NN$,
\beq\label{estimvek}
\left\|v_\ep^k-\int_0^{t}w\,ds\right\|^2_{L^\infty(0,T;H^1_0(\Om))^N}+\big\|\partial_tv_\ep^k-w\big\|^2_{L^\infty(0,T;L^2(\Om))^N}+k\,\big\|\partial_tv_\ep^k-w\big\|^2_{L^2(Q)}\leq C_w.
\eeq
\end{Lem}
\noindent
{\bf Proof of Theorem \ref{thmg}.}
Let $\{w^n,n\in\NN\}$ be a subset of $C^\infty_c(Q)^N$ which is dense in $L^1(0,T;L^2(\Om))^N$.
Let $u_\ep^{k,n}$, $k,n\in\NN$, be the solution to
\beq\label{pbawn}
\left\{\ba{l}\dis \rho\,\partial^2_{tt}u^{k,n}_\ep-{\rm Div}_x\big(Ae(u^{k,n}_\ep)\big)+(F_\ep+G_\ep) \partial_t u^{k,n}_\ep+k\big(\partial_t u^{k,n}_\ep-w^n\big)=0\ \hbox{ in }Q
\\ \ecart
\dis u^{k,n}_\ep=0\ \hbox{ on }(0,T)\times\partial\Om
\\ \ecart
u^{k,n}_\ep(0,.)=0,\ \partial_t u^{k,n}_\ep(0,.)=0\ \mbox{ in }\Om.
\ea\right.
\eeq
By virtue of Theorem \ref{thmS3} and using a diagonal extraction procedure, there exists a subsequence of $\ep$, still denoted by $\ep$, such that
the following convergences hold for any $k,n\in\NN$,
\beq\label{convuknep}
\left\{\ba{l}
u^{k,n}_\ep\stackrel{\ast}\rightharpoonup u^{k,n}\ \hbox{ in }L^\infty(0,T;H^1_0(\Om))^N\cap W^{1,\infty}(0,T;L^2(\Om))^N,
\\ \ecart
G_\ep\partial_t u^{k,n}_\ep\stackrel{\ast}\rightharpoonup g^{k,n}\ \hbox{ in }L^\infty(0,T;L^2(\Om))^N,
\ea\right.
\eeq
where $u^{k,n}$ is a solution to
\beq\label{pbukn}
\left\{\ba{l}
\dis (\rho I_N+M)\partial^2_{tt}u^{k,n}-{\rm Div}_x\big(Ae(u^{k,n})\big)+k(\partial_t u^{k,n}-w^n)+g^{k,n}=f\ \hbox{ in }Q
\\ \ecart
\dis u^{k,n}=0\ \hbox{ on }(0,T)\times\partial\Om
\\ \ecart
\dis u^{k,n}(0,.)=0,\ \partial_t u^{k,n}(0,.)=0\ \mbox{ in }\Om.
\ea\right.
\eeq
Fix $n\in\NN$. By the first convergence of \eqref{convuknep} and the estimate \eqref{estimvek} with $v^k_\ep=u^{k,n}_\ep$ and $w=w^n$, we have
\beq\label{convdtukn}
\partial_t u^{k,n}\mathop{\longrightarrow}_{k\to\infty} w^n\ \mbox{ in }L^2(Q)^N.
\eeq
Moreover, since the initial conditions of \eqref{pbawn} are clearly well-prepared in the sense \eqref{condbp}, by estimate \eqref{Estifg2} with $g=g^{k,n}$, we have for any $\bar x\in \bar\Om$, any $S\in (0,T)$ and a.e. $s\in (0,S)$,
\[
\int_{B(\bar x,S,s)}|g^{k,n}|^2dx\leq C\left(\int_0^s\left(\int_{B(\bar x,S,t)}|\partial_t u^{k,n}|^2dx\right)^{1\over 2}dt\right)^2,
\]
where the constant $C$ only depends on $\dis\sup_{\ep>0}\|G_\ep\|_{L^\infty(Q)^{N\times N}}$.
This combined with \eqref{convdtukn} yields
\[
\limsup_{k\to\infty}\,\int_{B(\bar x,S,s)}|g^{k,n}|^2dx\leq C\left(\int_0^s\left(\int_{B(\bar x,S,t)}|w^n|^2dx\right)^{1\over 2}dt\right)^2.
\]
Hence, using a diagonal extraction argument, there exist a subsequence of $k$, still denoted by $k$, such that for any $n\in\NN$,
\beq\label{convgkn}
g^{k,n}\stackrel{\ast}\rightharpoonup g^n\ \hbox{ in }L^\infty(0,T;L^2(\Om)))^N,
\eeq
which implies that for any $\bar x\in \bar\Om$, any $S\in (0,T)$ and a.e. $s\in (0,S)$,
\beq\label{estgn}
\int_{B(\bar x,S,s)}|g^n|^2dx\leq C\left(\int_0^s\left(\int_{B(\bar x,S,t)}|w^n|^2dx\right)^{1\over 2}dt\right)^2.
\eeq
Then, for any $w\in L^1(0,T;L^2(\Om))^N$ and any subsequence $w^{p_n}$ which converges strongly to $w$, we define the function $\G w$ by
\beq\label{defG}
g^{p_n}\stackrel{*}{\rightharpoonup} \G w\ \hbox{ in }L^\infty(0,T;L^2(\Om)))^N.
\eeq
This definition is independent of the strongly convergent subsequence $w^{p_n}$ due to the linearity of \eqref{pbawn} combined with estimate \eqref{estgn}.
By the linearity of problem \eqref{pbawn} the operator $\G$ is linear.
Moreover, using the lower semicontinuity of the $L^2(\Om)^N$-norm in \eqref{estgn} we deduce that $\G$ satisfies estimate \eqref{EstimateG1}. Estimate \eqref{EstimateG1b} is a simple consequence of \eqref{Estifg2b}  in the absence of measure $\mu^0$.
\par
Note that the definition of $\G$ is based on the subsequence $\ep$ satisfying convergences \eqref{convuknep} for any $k,n\in\NN$.
\par
Now let us prove estimate \eqref{EstimateG2}. Let $u_\ep$ be the solution to problem \eqref{pb1} and consider a subsequence $\ep'$ of $\ep$  such that $u_{\ep'}$ satisfies the results of Theorem~\ref{thmS3}. Also consider a sequence $w^{p_n}$ which strongly converges to $\partial_t u$ in $L^2(Q)^N$. Applying the estimate \eqref{Estifg2} with the sequence $u_{\ep'}-u^{k,p_n}_{\ep'}$ for $k,n\in\NN$, we get that for any $\bar x\in \bar\Om$, any $S\in (0,T)$ and a.e. $s\in (0,S)$,
\beq\label{estggkn}
\int_{B(\bar x,S,s)}|g-g^{k,p_n}|^2dx\leq C\mu^0\big( \bar B(\bar x,S,0)\big)+C\left(\int_0^s\left(\int_{B(\bar x,S,t)}|\partial_tu-\partial_t u^{k,p_n}|^2dx\right)^{1\over 2}dt\right)^2.
\eeq
where the measure $\mu^0$ is defined by \eqref{wepococ} with the sequence $u_{\ep'}$ but independently of $u^{k,p_n}_{\ep'}$.
Therefore, passing successively to the limit $k\to\infty$ with convergences \eqref{convgkn} and \eqref{convdtukn}, then to the limit $n\to\infty$ with convergences \eqref{defG} and $w^{p_n}\to\partial_t u$, we obtain the desired estimate \eqref{EstimateG2}.
This concludes the proof of Theorem~\eqref{thmg}.
\cqfd
\par\medskip\noindent
{\bf Proof of Corollary \ref{Thbcin}.}
Consider a subsequence of $\ep$ such that \eqref{convuep}, \eqref{Defig} and \eqref{Defzet} hold. Since \eqref{condbp} is satisfied, the function $\zeta$ defined by \eqref{Defzet} agrees with $u^1$ and the measure $\mu^0$ defined by \eqref{wepococ} vanishes. By \eqref{EstimateG2} we get that $g=\G (\partial_tu)$, and thus \eqref{pb2} proves that $u$ is a solution to \eqref{pb2Cibp}. Estimates \eqref{estimacot2}, \eqref{EstimateG1b} and Gronwall's Lemma imply the uniqueness of a solution to \eqref{pb2Cibp}.
Hence, it is not necessary to extract a new subsequence to get the convergence of $u_\ep$.
\cqfd
\par\medskip\noindent
{\bf Proof of Corollary \ref{Thbcinf}.} Consider the subsequence of $\ep$ given by Theorem \ref{thmg} and a dense countable set $\{\varphi^1_k\}$ of $L^2(\Om)^N$ contained in $C^\infty_c(\Om)^N$. By Theorem~\ref{thmS3}, Theorem~\ref{thmg} and \eqref{wepococ} we can use a diagonal argument to deduce the existence of a subsequence of $\ep$ and a linear operator $\F:{\rm Span}(\{\varphi^1_k\})\to L^\infty(0,T;L^2(\Om))$ such that for every $\varphi^1\in {\rm Span}(\{\varphi^1_k\})$ the solution $v_\ep$ of 
\beq\label{DCocf1}
\left\{\ba{l}\dis \rho\,\partial^2_{tt}v_\ep-{\rm Div}_x\big(Ae(v_\ep)\big)+B_\ep \partial_t v_\ep=0\ \hbox{ in }Q
\\ \ecart
\dis v_\ep=0\ \hbox{ on }(0,T)\times\partial\Om
\\ \ecart
\dis v_\ep(0,.)=0,\ \partial_t v_\ep(0,.)=\varphi^1\ \mbox{ in }\Om,
\ea\right.
\eeq
converges weakly-$\ast$ in $L^\infty(0,T;H^1_0(\Om))^N\cap W^{1,\infty}(0,T;L^2(\Om))^N$ to a function $v$ solution to
\beq\label{DCocf3}
\left\{\ba{l}\dis \rho\,\partial^2_{tt}v-{\rm Div}_x\big(Ae(v)\big)+B \partial_tv+\G(\partial_tv)=\F(\varphi^1)\ \hbox{ in }Q
\\ \ecart
\dis v=0\ \hbox{ on }(0,T)\times\partial\Om
\\ \ecart
\dis v(0,.)=0,\ \partial_t v(0,.)=\rho(\rho I_N+M)^{-1}\varphi^1\ \mbox{ in }\Om,
\ea\right.
\eeq
where $\F (\varphi^1)$ satisfies
\beq\label{DCocf2}
G_\ep\partial_t v_\ep\stackrel{\ast}\rightharpoonup \G(\partial_tv)-\F(\varphi^1)\ \hbox{ in }L^\infty(0,T;L^2(\Om))^N.
\eeq
Moreover,  by \eqref{wepococ} and estimate \eqref{EstimateG2} we have for any $S\in (0,T)$ and a.e. $s\in (0,S)$,
\beq\label{DCocf4}
\int_{B(\bar x,S,s)}\big|\F(\varphi^1)\big|^2dx\leq C\rho\int_{B(\bar x,S,0)}(\rho I_N+M)^{-1}M\varphi^1\cdot\varphi^1\,dx.
\eeq
This allows us to extend $\F$ to a continuous linear operator in $L^2(\Om)$ which satisfies \eqref{estimopP}.
\par
Assume now $u^0\in H^1_0(\Om)^N$, $u^1\in L^2(\Om)^N$ and define $u_\ep$ as the solution to \eqref{pb1} with $u^0_\ep=u^0$, $u^1_\ep=u^1$. Applying  Theorem~\ref{thmS3} and Theorem~\ref{thmg}, we can extract a subsequence of $\ep$ satisfying \eqref{convuep} and \eqref{Defig}, where $u$ is a solution to \eqref{pb2} with $\zeta=0$. Also applying Theorem \ref{thmg} to the sequence $u_\ep-v_\ep$, where $v_\ep$ is the solution to \eqref{DCocf1} for some $\varphi^1\in {\rm Span}(\{\varphi^1_k\})$, and recalling the definition \eqref{DCocf2} of $\F(\varphi^1)$, we have  for any $S\in (0,T)$ and a.e. $s\in (0,S)$,
\[
\int_{B(\bar x,S,s)}\big|g-\G(\partial_t u)+\F(\varphi^1)\big|^2dx\leq C\rho\int_{B(\bar x,S,0)}(\rho I_N+M)^{-1}M(u^1-\varphi^1)\cdot(u^1-\varphi^1) \,dx,
\]
which by the arbitrariness of $\varphi^1$ shows that
\beq\label{DCocf5}
g=\G(\partial_tu)-\F(u^1),
\eeq
and thus that $u$ is a solution to \eqref{pb2CipP}.
\par
The uniqueness of a solution to \eqref{pb2CipP} just follows by the uniqueness of a solution to \eqref{pb2Cibp} proved above, where $f$ is now replaced by $f+\F(u_1)$. This shows that it is not necessary to extract a new subsequence.
\par
Finally, estimate \eqref{positp} is a consequence of \eqref{Estifg2b} and \eqref{DCocf5}. \cqfd
\subsection{A general representation result}
The operator $\G$ defined by \eqref{defG} admits the following representation  which shows explicitly that $\G(\partial_t u)$ is a nonlocal operator with respect to the velocity in Theorem~\ref{thmg}.
\begin{Thm}\label{thmrep}
Under the assumptions of Theorem \ref{thmg} there exists a matrix-valued measure $\Lambda\in \M(\bar{Q};L^2(Q))^{N\times N}$ which is absolutely continuous with respect to the Lebesgue measure, such that $L^2(Q)^N\subset L^1(Q;d\Lambda)$ and such that the operator $\G$ defined by \eqref{defG} satisfies the representation formula
\beq\label{repG}
\G(w)(t,x)=\int_{Q}\,d\Lambda(s,y)\,w(s,y)\ \mbox{ a.e. in } Q,\quad\forall\,w\in L^2(Q)^N.
\eeq
Moreover, we have
\beq\label{suppLa}
\Lambda(B)=0\ \hbox{ a.e. in }\big\{(\bar x,S)\in \bar\Om\times (0,T): \ |B\cap \C(\bar x,S)|=0\big\},\quad\forall\, B\subset Q,\hbox{ measurable.}
\eeq
\end{Thm}
\par
Theorem \ref{thmrep} is based on the following representation result with Remark~\ref{remLambda} below.
\begin{Pro}\label{ProRep}
Let $X$ be a reflexive Banach space and let $(\om,\Sigma, \mu)$ be a finite measurable space. Then, for any linear continuous operator $\T: L^p(\om;d\mu)^N\to X$ with $1\leq p<\infty$, there exists a vector-valued measure
$\Lambda\in\M(\om;X)^N$, $\Lambda=(\Lambda_1,\ldots,\Lambda_N)$ which is absolutely continuous with  respect to $\mu$ such that 
\beq\label{intu}
L^p(\om;d\mu)^N\subset L^1(\om;\Lambda),\quad \T u=\int_\om d\Lambda(y)\,u(y),\quad \forall\, u\in L^p(\om;d\mu)^N,
\eeq
and
\beq\label{desNpr}
\sup_{E\in\Sigma}\|\Lambda_j(E)\|_X\leq |\om|^{1\over p}\|\T^\ast_j\|_{\L(X';L^{p'}(\om))}\leq 4\sup_{E\in\Sigma}\|\Lambda_j(E)\|_X,\quad \forall\,j\in\{1,\dots,N\},
\eeq
where we have
\[
\|\T_l^\ast\|\leq \|\T\|\leq\left(\sum_{j=1}^N\|\T_j^\ast\|^{p'}\right)^{1\over p'},\quad \forall\, l\in\{1,\dots,N\}.
\]
\end{Pro}
\begin{Rem}\label{remLambda}
We are mainly interested in the case where  $X=L^q(\varpi;d\nu)^M$ with $1<q<\infty$
In this case, $\Lambda=(\Lambda_1,\ldots,\Lambda_N)$ is replaced by a matrix-valued measure
\[
\Lambda_j=(\Lambda_{1j},\ldots, \Lambda_{Mj})\in \M(\om;L^q(\varpi;d\nu))^M,\quad \forall\,j\in\{1,\dots,M\}.
\]
Thus, $\Lambda$ belongs to $\M(\om;L^q(\varpi;d\nu))(\om;d\mu)^{M\times N}$ and \eqref{intu} can be written as
\beq\label{intuM}
L^p(\om;d\mu)^N\subset L^1(\Om;\Lambda),\quad \T u=\int_\om d\Lambda(y)\,u(y),\quad \forall\, u\in L^p(\om;d\mu)^N,
\eeq
where 
\[
\left(\int_\om  d\Lambda(y)u(y)\right)_j=\sum_{k=1}^N\int_\om d\Lambda_{jk}(y)\,u_k(y),\quad \forall\,j\in\{1,\dots,M\}.
\]
\par
Observe that for any set $E\subset \Sigma$, $\Lambda(E)$ is a function in $L^q(\varpi;d\nu)^{M\times N}$, then the $M\times N$ matrix $\Lambda(E)(x)$ is defined $\nu$-a.e. $x\in \varpi$. If we assume that
\beq\label{HypKer}
\exists\,\N\subset \varpi,\ \mu(\N)=0,\ \hbox{such that }
\left\{\ba{l}\dis\hbox{the function }
E\in\Sigma\mapsto \Lambda(E)(x)
\\ \ecart
\dis \hbox{is well defined for any }x\in \varpi\setminus \N\hbox{and defines a measure,}\ea\right.
\eeq
then, denoting $\Lambda(x,E)=\Lambda(E)(x)$, formula \eqref{intuM} can be written as the kernel representation formula
\beq\label{intuMx}
(\T u)(x)=\int_\om d\Lambda(x,y)\,u(y),\quad \forall\, u\in L^p(\om;d\mu)^N.
\eeq
However, it is not clear than assumption \eqref{HypKer} holds true in general.
Furthermore, even if formula \eqref{intuMx} holds, $\Lambda(x,.)$ is not in general absolutely continuous with respect to $\mu$, {\em i.e.} $\Lambda(x,.)$ is not a function but just a measure. As a simple example,
consider $L^q(\varpi;d\nu)^M= L^p(\om;d\mu)^N$ and $\T$ as the identity operator, then the measure $\Lambda$ is given by
\[
\Lambda(B)=1_B\,I_N,\ \mbox{ for }B\in \Sigma.
\]
In this case \eqref{intuMx} is satisfied with
\[
\Lambda(x,y)=\delta_x(y)\,I_N.
\]
\end{Rem}
\noindent
{\bf Proof of Theorem \ref{thmrep}.}
First note that
\[
\G: L^2(Q)^N\hookrightarrow L^1(0,T;L^2(\Om))^N\longrightarrow L^\infty(0,T;L^2(\Om))^N\hookrightarrow L^2(Q)^N,
\]
where the two embedding are continuous. Moreover, by the Cauchy-Schwarz inequality and estimate \eqref{EstimateG1} we get that for any $w\in L^2(Q)^N$,
\beq\label{estGcone}
\int_{\C(\bar{x},S)}\big|\G w\big|^2dt\,dx\leq {1\over 2}\,C S^2\int_{\C(\bar{x},S)}|w|^2dt\,dx,\quad\forall\,(S,\bar{x})\in Q.
\eeq
which implies in particular the continuity of the linear operator $\G$ from $L^2(Q)^N$ into $L^2(Q)^N$.
\par
Therefore, applying Proposition \ref{ProRep} and Remark \ref{remLambda} with $X=L^2(Q)^N$, $\om=Q$, $\mu$ the Lebesgue measure on $Q$ and $p=2$, there exists a matrix-valued measure $\Lambda\in \M(\bar{Q};L^2(Q))^{N\times N}$ which is absolutely continuous with respect to the Lebesgue measure, such that $\G$ satisfies the representation formula \eqref{repG}.
\par
Moreover, applying \eqref{EstimateG1} we get \eqref{suppLa}.
\cqfd
\par\bigskip\noindent
{\bf Proof of Proposition \ref{ProRep}.}
Denoting $i_{p',1}$ the continuous embedding from $L^{p'}(\om;d\mu)^N$ into $L^1(\om;d\mu)^N$, we apply  Theorem 8.1 in \cite{DuSc} to the $N$ components of the operator $i_{p'1}\circ \T^\ast$ in ${\L}(X';L^1(\om;d\mu))^N$. Taking into account that $X'$ is reflexive and then the unit ball is weakly compact, we deduce that there exists a vector-valued measure $\Lambda=(\Lambda_1,\ldots,\Lambda_N)\in \M(\om;X)^N$, which is absolutely continuous with respect to $\mu$, such that  for any $\zeta'\in X'$ and any $j\in \{1,\ldots,N\}$, the measure $E\in\Sigma\mapsto\big\langle \zeta',\Lambda_j(E)\big\rangle_{X',X}\in\RR$ satisfies
\[
\T^\ast_j(\zeta')={d\over d\mu}\,\big\langle \zeta',\Lambda_j(.)\big\rangle_{X',X},
\]
or equivalently
\[
\int_E (\T^\ast_j\zeta')(x)\,d\mu(x)= \big\langle \zeta' ,\Lambda_j(E)\big\rangle_{X',X},\quad \forall\, E\in \Sigma.
\]
Therefore, for any step function 
\[
u=\sum_{l=1}^m\lambda_l\,1_{E_l},\quad \lambda_1,\ldots,\lambda_m\in \RR^N,\ E_1,\ldots,E_m\in \Sigma,
\]
and any $\zeta'\in X'$, we have
\[
\ba{l}
\dis \langle\zeta',\T u\rangle_{X',X}=\big\langle\T^\ast\zeta',u\big\rangle_{L^{p'}(\om)^N,L^{p}(\om)^N}=\sum_{j=1}^N\int_\om (\T^\ast_j\zeta')(x)u_j(x)\,d\mu(x)
\\ \ecart
\dis =\sum_{j=1}^N\sum_{l=1}^m \lambda_{lj}\int_{E_l}(\T^\ast_j\zeta')(x)\,d\mu(x)=\sum_{j=1}^N\sum_{l=1}^m \lambda_{lj}\big\langle \zeta',\Lambda_j(E_l)\big\rangle_{X',X}
\\ \ecart
\dis =\sum_{l=1}^m\big\langle \zeta',\lambda_l\cdot \Lambda(E_l)\big\rangle_{X',X}=
\left\langle\zeta',\int_\om d\Lambda(y)\,u(y)\right\rangle.
\ea
\]
This shows that for any step function $u$,
\[
\T u=\int_\om d\Lambda(y)\,u(y).
\]
Now, using that $\T$ is a continuous operator from $L^p(\om;d\mu))^N$ into $X$, we conclude to \eqref{desNpr}.
\cqfd
\subsection{Proof of the lemmas}
\noindent
{\bf Proof of Lemma \ref{LemDeriv}.}
Using a translation we can always assume that $x_0=0$.
\par
First assume $B_{R(t)}\subset\Om$, for any $t\in [0,T]$. In this case, using the change of variables
$y=x/R(t)$, we have
\[
\Phi(t)=R(t)^N\int_{B_1}u(t,R(t)y)\,dy.
\]
Then, we have
\[
\ba{ll}\dis\Phi'(t) &\dis=NR(t)^{N-1}R'(t)\int_{B_1}u(t,R(t)y)\,dy+R(t)^N\int_{B_1}\partial_tu(t,R(t)y)\,dy
\\ \ecart
& \dis + R(t)^NR'(t)\int_{B_1}\nabla_xu(t,R(t)y)\cdot y\,dy =R(t)^N\int_{B_1}\partial_tu(t,R(t)y)\,dy
\\ \ecart
& \dis + R(t)^{N-1}R'(t)\int_{B_1}\big(Nu(t,R(t)y)+\nabla_y\big[u(t,R(t)y)\big]\cdot y\big)dy
\\ \ecart &
\dis =R(t)^N\int_{B_1}\partial_tu(t,R(t)y)\,dy+R(t)^{N-1}R'(t)\int_{B_1}{\rm div}_y\big[u(t,R(t)y)y\big]dy
\\ \ecart
& \dis =R(t)^N\int_{B_1}\partial_tu(t,R(t)y)\,dy+R(t)^{N-1}R'(t)\int_{\partial B_1}u(t,R(t)y)\,dy
\\ \ecart
& \dis =\int_{B_{R(t)}}\partial_tu(t,x)\,dx+R'(t)\int_{\partial B_{R(t)}}u(t,x)\,dx,
\ea
\]
which proves the result. \par
In the general case, by the regularity of $\Om$ we can always assume that
\[
u\in W^{1,1}(0,T;L^1(\RR^N))\cap L^1(0,T;W^{1,1}(\RR^N)).
\]
For  $\lambda\in C^\infty_c(\RR^N)$ with
\[
\int_{\RR^N}\lambda(x)\,dx=1,
\]
and $\ep>0$, define $\zeta_\ep\in C^\infty(\RR^N)$ by
\[
\zeta_\ep(x):={1\over \ep^N}\int_{\Om}\lambda\left({x-y\over \ep}\right)\,dy\quad\forall\, x\in\RR^N.
\]
Applying the above proved to the function $(t,x)\mapsto u(t,x)\zeta_\ep(x)$, we get
\beq\label{eqzetaep}
\ba{l}\dis \int_{B_{R(t_2)}}\hskip-10pt u(t_2,x)\,\zeta_\ep(x)\,dx-\int_{B_{R(t_1)}}\hskip-10pt u(t_1,x)\,\zeta_\ep(x)\,dx
\\ \ecart
\dis =\int_{t_1}^{t_2}\left(\int_{B_{R(t)}}\hskip-10pt \partial_tu(t,x)\,\zeta_\ep(x)\,dx+R'(t)\int_{\partial B_{R(t)}} \hskip-10pt  u(t,x)\,\zeta_\ep(x)\,ds(x)\right)dt.
\ea
\eeq
Moreover, using that $\zeta_\ep$ is bounded in $L^\infty(\RR^N)$ and 
\[
\zeta_\ep(x)\to\left\{\ba{ll} 1 &\hbox{ if }x\in \Om\\ \ecart\dis{1\over 2}&\hbox{ if }x\in \partial\Om
\\ \ecart
\dis 0 &\hbox{ if }x\in \RR^N\setminus\bar\Om,\ea\right.\quad \hbox{ when }\ep\to 0,
\]
we can pass to the limit in \eqref{eqzetaep} to deduce
\beq\label{ecu1deiv}
\ba{l}\dis \int_{B_{R(t_2)}\cap\Om}\hskip-10pt u(t_2,x) dx-\int_{B_{R(t_1)}\cap\Om}\hskip-10pt u(t_1,x) dx
\\ \ecart
\dis =\int_{t_1}^{t_2}\left(\int_{B_{R(t)}\cap\Om}\hskip-1.em \partial_tu(t,x)\,dx+R'(t)\int_{\partial B_{R(t)}\cap\Om} u(t,x)\,ds(x)+{1\over 2}R'(t)\int_{\partial B_{R(t)}\cap\partial\Om} u(t,x)\,ds(x)\right)dt.
\ea
\eeq
\par
Now, consider  $(\bar t,\bar x)\in (0,T)\times\partial\Om$ such that
\[
R'(\bar t)\not =0,\quad \bar x\in\partial B_{R(\bar t)}\cap \partial\Om.
\]
Since $\Om$ is $C^1$-regular, there exists a ball $B(\bar x,\delta_{\bar x})$, an open set $O\subset \RR^{N-1}$ with $0\in O$ and a function $\phi=\phi(\zeta)\in C^1(O;\RR^N)$ such that 
\[
\phi(0)=\bar x,\quad \phi\hbox{ is injective in }O,\quad \hbox{Rank}(D\phi)(\zeta)=N-1,\ \forall\, \zeta\in O,\quad \partial\Om\cap B(\bar x,\delta_{\bar x})=\phi(O).
\]
Since $R'(\bar t)\not=0$, applying the implicit function theorem to the function
\[
(t,\zeta)\in (0,T)\times O\mapsto |\phi(\zeta)|-R(t),
\]
$\delta_{\bar x}$ and $O$ can be chosen small enough to ensure the existence of $\ep>0$ and a function $\psi$ in $C^1(O;(\bar t-\ep,\bar t+\ep))$ such that 
\beq\label{ConTim}
\psi(0)=\bar t,\quad R(\psi(\zeta))=|\phi(\zeta)|,\ \forall\, \zeta\in O,\quad \left\{\ba{l}\dis t\in (\bar t-\ep,\bar t +\ep),\ \zeta\in O\\ \ecart\dis  R(t)=|\phi(\zeta)|\
\ea\right.\Rightarrow\; t=\psi(\zeta).
\eeq
Therefore, we have
\[
\big\{(t,x):\ t\in (\bar t-\ep,\bar t+\ep):x\in B(\bar x,\delta_{\bar x})\cap \partial B_{R(t)}\cap\partial\Om\}=\big\{\big(\psi(\zeta),\zeta\big):\zeta\in O\big\},
\]
which  thus has null $N$-dimensional measure.
Since for any integer $n\geq 1$, the set (we can assume that $R$ is defined in an interval larger than $[0,T]$)
\[
\big\{(\bar t,\bar x)\in [0,T]\times\partial\Om:\bar x\in \partial B_{R(t)},\,|R'(t)|\geq 1/n\big\}
\]
is a compact, we deduce that this set has zero measure.
Hence,
\[
\big\{(\bar t,\bar x)\in [0,T]\times\partial\Om:\bar x\in \partial B_{R(t)},\,R'(t)\not =0\big\}
\]
also has null $N$-dimensional measure.
This shows that the last term in \eqref{ecu1deiv} is zero for a.e. $t\in (0,T)$.
Therefore, the derivative formula \eqref{derivinv} holds.
\cqfd
\par\medskip\noindent
{\bf Proof of Lemma \ref{lemsecef}.}  It is enough to  consider the case where $\varphi\in C^1(\bar Q)$ and $\nu\in C^1(\bar Q)$.
For any integer $n\geq 1$ and any $k\in \{0,\ldots,n-1\}$, set
\beq\label{moyti}
\bar u_\ep^{n,k}(x):={n\over T}\int_{{k\over n}T}^{{k+1\over n}T}u_\ep(t,x) dt,\quad \bar v_\ep^{n,k}(x):={n\over T}\left(u_\ep\left({k+1\over n}T,x\right)-u_\ep\left({k\over n}T,x\right)\right),
\eeq
\beq\label{moytil}
\bar u^{n,k}:={n\over T}\int_{{k\over n}T}^{{k+1\over n}T}u(t,x) dt,\quad \bar v^{n,k}:={n\over T}\left(u\left({k+1\over n}T,\cdot\right)-u\left({k\over n}T,\cdot\right)\right),
\eeq
\beq\label{moyph}
\bar \varphi^{n,k}(x):={n\over T}\int_{{k\over n}T}^{{k+1\over n}T}\varphi(t,x) dt,\quad  \bar \nu^{n,k}(x):={n\over T}\int_{{k\over n}T}^{{k+1\over n}T}\nu (t,x)dt.
\eeq
\par
Taking into account that for any $\xi\in \RR^N$ with $|\xi|\leq 1$,  the function ${\mathscr Q}_\xi$ defined as
\[
{\mathscr Q}_\xi:(v,V)\in\RR^N\times\RR^{N\times N}_s\;\longmapsto\; {c\over 2}\rho|v|^2+{c\over 2}AV:V-AV:v\odot \xi\in\RR
\]
is convex and \eqref{corresCDT}, we have 
\[
\ba{l}\dis \int_0^{T}\int_\Om{\mathscr Q}_{\nu}\big(\partial_tu_\ep,e(u_\ep)\big)\varphi\,dx\,dt=\sum_{j=0}^{n-1}\int_{{k\over n}T}^{{k+1\over n}T}\int_\Om{\mathscr Q}_{\bar\nu^{n,k}}\big(\partial_tu_\ep,e(u_\ep)\big)\bar \varphi^{n,k}dx\,dt-{C\over n}
\\ \ecart
\dis\geq{T\over n}\,\sum_{k=0}^{n-1}\int_\Om{\mathscr Q}_{\bar\nu^{n,k}}\big(\bar v_\ep^{n,k},e(\bar u_\ep^{n,k})\big)\bar\varphi^{n,k}dx\,dt-{C\over n}
\\ \ecart
\dis={T\over n}\,\sum_{k=0}^{n-1} \int_\Om{\mathscr Q}_{\bar\nu^{n,k}}\Big(\bar v_\ep^{n,k},e(\bar u^{n,k})+\sum_{j=1}^Ne(w_\ep^j)\bar v_{j,\ep}^{n,k}\Big)\bar\varphi^{n,k}dx+O_\ep-{C\over n}
\\ \ecart\dis
={T\over n}\,\sum_{k=0}^{n-1} \int_\Om{\mathscr Q}_{\bar\nu^{n,k}}\Big(\bar v^{n,k}+(\bar v_\ep^{n,k}-\bar v^{n,k}),e(\bar u^{n,k})+\sum_{j=1}^Ne(w_\ep^j)(\bar v_{j}^{n,k}+(\bar v_\ep^{n,k}-\bar v^{n,k}))\Big)\bar\varphi^{n,k}dx+O_\ep-{C\over n}.
\ea
\]
Using the weak convergence to zero of $e(w_\ep^j)$ in $L^p(\Om)^{N\times N}$ with $p>N$, Rellich-Kondrachov's compactness theorem for $\bar v_\ep^{n,k}-\bar v^{n,k}$, the non-negativity of the quadratic form ${\mathscr Q}_\nu$ and the definition \eqref{defM} of $M$, we have
\[
\ba{l}\dis \int_0^{t}\int_\Om{\mathscr Q}_{\nu}\big(\partial_tu_\ep,e(u_\ep)\big)\varphi\,dx\,dt\geq
{T\over n}\,\sum_{k=0}^{n-1} \int_\Om{\mathscr Q}_{\bar\nu^{n,k}}\Big(\bar v^{n,k},e(\bar u^{n,k})+\sum_{j=1}^Ne(w_\ep^j)\bar v_{j}^{n,k}\Big)\bar\varphi^{n,k}dx+O_\ep-{C\over n}
\\ \ecart
\dis ={T\over n}\,\sum_{k=0}^{n-1} \int_\Om\Big({\mathscr Q}_{\bar\nu^{n,k}}\big(\bar v^{n,k},e(\bar u^{n,k})\big)+{c\over 2}M\bar v^{n,k}\cdot\bar v^{n,k}\Big)\bar\varphi^{n,k}dx+O_\ep-{C\over n}.
\ea
\]
Therefore, we have just proved
\[
\liminf_{\ep\to 0}\int_0^{T}\int_\Om{\mathscr Q}_{\nu}\big(\partial_tu_\ep,e(u_\ep)\big)\varphi\,dx\,dt\geq
{T\over n}\,\sum_{k=0}^{n-1}\int_\Om\Big({\mathscr Q}_{\bar\nu^{n,k}}\big(\bar v^{n,k},e(\bar u^{n,k})\big)+{c\over 2}M\bar v^{n,k}\cdot\bar v^{n,k}\Big)\bar\varphi^{n,k}dx-{C\over n}.
\]
Passing to the limit as $n$ tends to infinity, we finally obtain \eqref{regdtu} and \eqref{semlf}.
\par
Finally, let us prove \eqref{liuews}.
The sequence $Ae(u_\ep):e(w_\ep^j)$ is bounded in $L^\infty\big(0,T;L^{{2p\over p+2}}(\Om)\big)$.
Hence, it is enough that convergence \eqref{liuews} holds in the distributions sense in $Q$.
\par
Let $\ph\in C^\infty_c(Q)$.
With notations \eqref{moyti} and \eqref{moytil} we have for any integer $n\geq 1$,
\[
\ba{l}
\dis \int_Q Ae(u_\ep):e(w_\ep^j)\,\ph\,dx\,dt=\sum_{k=0}^{n-1}\int_{{k\over n}T}^{{k+1\over n}T}\int_\Om Ae(u_\ep):e(w_\ep^j)\,\bar\varphi^{n,k}\,dx+O\left({1\over n}\right)
\\ \ecart
\dis ={T\over n}\,\sum_{k=0}^{n-1}\int_\Om Ae(\bar u_\ep^{n,k}):e(w_\ep^j)\,\bar\varphi^{n,k}\,dx+O\left({1\over n}\right).
\ea
\]
Then, again using convergence \eqref{corresCDT}, the weak convergence to zero of $e(w_\ep^j)$ in $L^p(\Om)^{N\times N}$ and the definition \eqref{defM} of $M$, we get
\[
\ba{l}
\dis \int_Q Ae(u_\ep):e(w_\ep^j)\,\ph\,dx\,dt={T\over n}\,\sum_{k=0}^{n-1}\int_\Om Ae\Big(\bar u^{n,k}+\sum_{i=1}^Ne(w_\ep^i)\bar v_i^{n,k}\Big):e(w_\ep^j)\,\bar\varphi^{n,k}\,dx+O_\ep+O\left({1\over n}\right).
\\ \ecart
\dis ={T\over n}\,\sum_{k=0}^{n-1}\int_\Om M\bar v^{n,k}\cdot e_j\,\bar\varphi^{n,k}\,dx\,dt+O_\ep+O\left({1\over n}\right).
\ea
\]
Passing successively to the limit as $\ep$ tends to zero for a fixed $n$, and to the limit as $n$ tends to infinity, we obtain that
\[
\lim_{\ep\to 0}\,\int_Q Ae(u_\ep):e(w_\ep^j)\,\ph\,dx\,dt=\int_Q M\partial_t u\cdot e_j\,\ph\,dx\,dt,
\]
which concludes the proof of Lemma~\ref{lemsecef}.
\cqfd
\par\medskip\noindent
{\bf Proof of Lemma \ref{lemthmg}.} Set
\[
v(t,.)=\int_0^{t} w(s,.)\,ds,\ \mbox{ for }t\in [0,T].
\]
Putting $\partial_t v_\ep^k-w$ as test function in \eqref{pbaw} and using that  $F_\ep$, $G_\ep$ are skew-symmetric, we get
\beq\label{estEpk1}
\ba{l}
\dis{1\over 2}{d\over dt}\into\Big(\rho\big|\partial_t v_\ep^k-w\big|^2+Ae(v_\ep^k-v):e(v_\ep^k-v)\Big)\,dx
+k\into \big|\partial_t v_\ep^k-w\big|^2dx
\\ \ecart
\dis ={d\over dt}\into F_\ep (v_\ep^k-v)\cdot w\,dx-\into F_\ep (v_\ep^k-v)\cdot\partial_tw\,dx-\into G_\ep w\cdot\big(\partial_t v_\ep^k-w\big)\,dx
\\ \ecart
\dis -\into \rho\,\partial_t w\cdot(\partial_tv_\ep^k-w)\,dx-\into Ae(v): e(\partial_t v_\ep^k-w)\,dx.
\ea\eeq
Setting
\[
E_\ep^k(t):={1\over 2}\,\into\Big(\rho\big|\partial_t v_\ep^k-w\big|^2+Ae(v_\ep^k-v):e(v_\ep^k-v)\Big)\,dx,
\]
\[
h_\ep^k(t):=\into F_\ep (v_\ep^k-v)\cdot w\,dx,
\]
equality \eqref{estEpk1} implies that
\beq\label{estEpk2}
{dE_\ep^k\over dt}+k\into \big|\partial_t v_\ep^k-w\big|^2dx\leq C_w\big(E_\ep^k+1\big)+{dh_\ep^k\over dt}.
\eeq
Applying Gronwall's lemma and noting that $E_\ep^k(0)=0$, we get
\[
E_\ep^k(t)\leq C_w+C_w\int_0^{t} E_\ep^k(s)\,ds,
\]
which again using Gronwall's lemma gives 
\beq\label{estEpk3}
E_\ep^k(t)\leq C_w,\quad \forall\,t\in [0,T].
\eeq 
This combined with \eqref{estEpk2} proves that
\[
k\int_Q \big|\partial_t v_\ep^k-w\big|^2dxdt\leq C_w.
\]
Finally, the former estimate and \eqref{estEpk3} yield the desired estimate \eqref{estimvek}.
\cqfd

\end{document}